\documentclass[11pt]{extarticle}

\usepackage[left=3cm,right=3cm,top=2.2cm,bottom=2.2cm]{geometry}
\usepackage{amsmath,amssymb,mathrsfs,stmaryrd}
\usepackage{amsthm}
\usepackage{mathtools,bm}
\usepackage{dsfont}
\usepackage[hidelinks]{hyperref}
\usepackage[nameinlink,capitalize]{cleveref}
\usepackage{xcolor}
\usepackage[shortlabels] {enumitem}
\usepackage{microtype}
\usepackage[abbrev,msc-links,alphabetic]{amsrefs}
\usepackage{todonotes}
\usepackage{authblk}
\usepackage{titlefoot}
\usepackage{float}
\usepackage{graphicx}

\theoremstyle{plain}
\newtheorem{theorem}{Theorem}[section]

\newtheorem{lemma}[theorem]{Lemma}

\newtheorem{assumption}[theorem]{Assumption}
\theoremstyle{definition}
\newtheorem{definition}[theorem]{Definition}

\newtheorem{problem}[theorem]{Problem}
\newtheorem{remark}[theorem]{Remark}

\colorlet{darkred}{red!90!black}

\newcommand{\E}{\ensuremath{{\mathbb E}}}
\renewcommand{\P}{\ensuremath{{\mathbb P}}}

\usepackage{graphicx}

\makeatletter
\newcommand*\bigcdot{{\mathpalette\bigcdot@{.5}}}
\newcommand*\bigcdot@[2]{\mathbin{\vcenter{\hbox{\scalebox{#2}{$\m@th#1\bullet$}}}}}
\makeatother

\newenvironment{proofThm}[1]
{
	\par\vspace{\baselineskip}\noindent
	\textit{Proof of Theorem} #1.
	\noindent
}
{
	\qed\ignorespacesafterend
}

\begin{document}
	
	\title{Approximation of Optimal Feedback Controls for Stochastic Reaction-Diffusion Equations}

	\author{Wilhelm Stannat}
	\author{Alexander Vogler}
	
	\nocite{*}
	
	\affil{\small Technische Universit\"at Berlin, Berlin, Germany}
	
	\maketitle

	\unmarkedfntext{\textit{Keywords and phrases ---} stochastic optimal control, reaction-diffusion equations, adjoint calculus, variational methods, artificial neural networks}
	
	\unmarkedfntext{\textit{Mail}: \textbullet$\,$ stannat@math.tu-berlin.de $\,$\textbullet$\,$ vogler@math.tu-berlin.de}
	
	\begin{abstract}
		In this paper we present a method to approximate optimal feedback controls for stochastic reaction-diffusion equations. We derive two approximation results providing the theoretical foundation of our approach and allowing for explicit error estimates. The approximation of optimal feedback controls by neural networks is discussed as an explicit application of our method. Finally we provide numerical examples to illustrate our findings.
	\end{abstract}
	
	
	\tableofcontents
	
	\section{Introduction}
	
	For fixed time horizon $T>0$ and bounded domain $\Lambda\subset \mathbb{R}$ we consider the controlled stochastic partial differential equation (SPDE)
	\begin{equation}\label{spde1}
		\begin{cases}
			du^{\mathfrak{g}}_t=[ A u^{\mathfrak{g}}_t+\mathcal{F}(u^{\mathfrak{g}}_t)+\mathfrak{g}_t]dt+ B dW_t, \quad t\in [0,T]\\
			u^{\mathfrak{g}}_0=u_0,
		\end{cases}
	\end{equation} 
	where $W:[0,T]\times \Omega\rightarrow \Xi$ denotes a cylindrical $Q$-Wiener process on a separable Hilbert space $\Xi$. The covariance operator $Q:\Xi\rightarrow \Xi$ is assumed to be linear, bounded, positive-definite and self-adjoint, $B\in L(\Xi,L^2(\Lambda))$ and $A:D(A)\subset L^2(\Lambda)\rightarrow L^2(\Lambda)$ is a densely defined self-adjoint, negative definite linear operator with domain $D(A)$ and compact inverse (for example the Dirichlet-Laplace). For a ONS  of eigenvectors $(e_n)_n$ in $L^2(\Lambda)$ with corresponding eigenvalues $(\lambda_n)_{n\geq 1}$ the domain of $A$ is characterized as
	\begin{align*}
		x\in D(A):=\{x=\sum_{n=1}^{\infty}x_ne_n\in L^2(\Lambda):\|x\|_2^2 :=\sum_{n=1}^{\infty}\lambda_nx_n^2<\infty\}.
	\end{align*}
	Our control problem can be formulated as follows:
	\begin{problem}\label{SCP}
		Minimize the cost functional 
		\begin{align*}
			J(\mathfrak{g}) =  J_1 (u^{\mathfrak{g}}) + J_2 ( \mathfrak{g} ),
		\end{align*}
		where
		\begin{align*}
			J_1(u)&:=\E\left[\int_{0}^{T}\int_{\Lambda}  l(t,x,u_t(x)) \mathrm{d}x dt + \int_{\Lambda}m(x,u_T(x))dx \right],\\
			J_2(\mathfrak{g})&:=\E\left[\int_{0}^{T}\| \mathfrak{g}_t\|^2_{L^2(\Lambda)} dt  \right],
		\end{align*}
		for running cost induced by $l:[0,T]\times \mathbb{R}\times \mathbb{R}\rightarrow \mathbb{R}_+$ and terminal cost induced by $m:\mathbb{R}\times \mathbb{R}\rightarrow \mathbb{R}_+$, subject to the SPDE \eqref{spde1}, over the set of admissible controls 
		\begin{align*}
			\mathbb{A}= L^2([0,T]\times \Omega,(\mathcal{F}_t)_{t\in [0,T]};\mathcal{U}),
		\end{align*}
		for a closed and convex subset $\mathcal{U}\subseteq L^2(\Lambda)$.
	\end{problem}
	
	The mathematical theory of these problems is by now well understood (see e.g. \cite{CDP18}, \cite{DM13}, \cite{FGS17a}, \cite{FZ20}, \cite{FO16}, \cite{FHT18}, \cite{LZ14}, \cite{LZ18}, \cite{SW20}, \cite{SW21}, \cite{SW22}), however the efficient numerical approximation of these control problems still faces serious difficulties due to the computational complexity of classical approaches that require either to approximate an infinite dimensional Hamilton-Jacobi-Bellman equation or a backward SPDE. 
	
	In the last couple of years there is a rising interest in more efficient methods for the numerical approximation of optimal control problems (\cite{BJ19}, \cite{CL22} \cite{DKK2021}, \cite{DMPV19}, ,\cite{EHJ17}, \cite{GKM18}, \cite{KK18}, \cite{NR21}, \cite{OSS22}, \cite{SW20}). In particular applications of machine learning algorithms to the approximation of backward stochastic differential equations and optimal control problems have drawn a lot of attention recently(\cite{EHJ17},\cite{CL21}, \cite{CL22}). However the literature in this topic for SPDE's is still very sparse. One main difficulty in the case of optimal control of SPDE's is that the solution to the state equation (\ref{spde1}) takes values in an infinite dimensional space. Assuming the optimal control $\hat{\mathfrak{g}}$ is of markovian feedback type, i.e. $\hat{\mathfrak{g}}_t=\hat{G}(t,u_t^{\hat{G}})$, the function $\hat{G}$ still depends through the solution on some infinite dimensional object, which makes it much more difficult to approximate using neural networks. 
	
	In our approach we consider finitely based approximations for the optimal feedback function $\hat{G}$, which are defined in terms of the projection of $\hat{G}$ onto some finite dimensional subspace. We show that it is indeed sufficient if one restricts the optimization to some space of admissible controls that can approximate the finitely based approximations of $\hat{G}$ in order to reach the optimal cost. Based on this observation we can construct spaces of admissible controls that allow for an efficient numerical approximation of the control problem. In the following we refer to those spaces as ansatz spaces. Our second main result provides explicit convergence rates under additional assumptions on the regularity of $\hat{G}$.
	
	The rest of this paper is organized as follows: In Section \ref{assumptions} we provide standing assumptions and state our main results. In Section \ref{thm1} we will give a prove for our first main theorem, which is  reminiscent of a universal approximation result. In Section \ref{thm2} we will prove our second main result, which provides convergence rates. The construction of ansatz spaces will be discussed in Section \ref{ansatz} and Section \ref{ex}. Finally Section \ref{numerics} is devoted to numerical examples.

	\section{Standing Assumptions and Main Results}\label{assumptions}
	
	\subsection{Standing Assumptions}
	
	In the following let $H:=L^2(\Lambda)$. Furthermore let $(e_n)_{n\ge 1}$ denote an orthonormal basis of $H$ consisting of eigenfunctions of $A$ with corresponding 
	eigenvalues $(-\lambda_n )_{n\ge 1}$ arranged in increasing order 
	$0< \lambda_1 \le \lambda_ 2\le \ldots $. 
	
	For $r\in \mathbb{R}$ we can then define the fractional operator 
	$(-A)^{\frac{r}{2}} : D((-A)^\frac{r}{2})\rightarrow L^2(\Lambda)$ by 
	\begin{align*}
		(-A)^{\frac{r}{2}}u := \sum_{n=1}^{\infty}\lambda_n^{\frac{r}{2}} 
		u_n e_n,
	\end{align*}
	for all 
	\begin{align*}
		u\in D((-A)^\frac{r}{2}) := \{u = \sum_{n=1}^{\infty} u_n e_n 
		\in H: \|u\|_r^2 := \sum_{n=1}^{\infty} \lambda^r_n
		u_n^2<\infty\}.
	\end{align*}
	To simplify notations let $H^r := D((-A)^\frac{r}{2})$ and note that 
	$\|u\|_r = \|(-A)^{\frac{r}{2}}u\|_{H}$ defines a norm on 
	$H^r$. By $L_2^0$ we denote the space of Hilbert-Schmidt operators 
	$\Phi:\Xi\rightarrow H$ with norm 
	\begin{align*}
		\|\Phi\|_{L_2^0}^2 := \sum_{n=1}^\infty\|\Phi \varphi_n\|_{H}^2,
	\end{align*} 
	for a basis $(\varphi_n)_{n\in \mathbb{N}}$ of $\Xi$, and for $\Phi$ taking values in $H^r$ we define 
	$\|\Phi\|_{L^0_{2,r}} := \|(-A)^{\frac{r}{2}}\Phi\|_{L_2^0}$.
	
	\begin{remark}
		If $A = \Delta$, where $\Delta$ denotes the Laplace operator with Dirichlet boundary conditions, it is well known that $H^1=H_0^1(\Lambda)$ and $H^2=H^2(\Lambda)\cap H_0^1(\Lambda)$, see e.g. \cite{Kru14}. In the case of Neumann boundary conditions one can consider the operator $(A-\alpha I)$ for some $\alpha >0$ and $\mathcal{F}(u)+\alpha u$ in the non-linearity instead.
	\end{remark}
	
	We impose the following assumptions on the coefficients of the controlled 
	spde \eqref{spde1}:
	
	\begin{assumption}\label{A1}
		\begin{enumerate}[H1.]
			\item We assume that $\cal F$ is of Nemytskii type, i.e. 
			$\mathcal{F} (u)(x) = f(u(x))$ for some $f:\mathbb{R}\rightarrow 
			\mathbb{R}$ continuously differentiable and Lipschitz continuous 
			\begin{align*}
				|f(u_1)-f(u_2)|\leq C|u_1-u_2| \quad \forall u_1 , u_2 \in \mathbb{R}, 
			\end{align*}
			for some constant $C>0$. 
			\item The dispersion operator $B\sqrt{Q}:\Xi\rightarrow H$ is Hilbert-Schmidt, i.e. $\|B\sqrt{Q}\|_{L_2^0}<\infty$. 
			\item For any $t\in [0,T]$ and $x\in \mathbb{R}$ the functions 
			$l(t,x,\cdot) :  \mathbb{R}\rightarrow 
			\mathbb{R}_+$ and $m(x,\cdot):\mathbb{R} \rightarrow 
			\mathbb{R}_+$ are differentiable.
			\item  For any $t\in [0,T]$ and $x\in \mathbb{R}$ the functions $l(t,x,\cdot)$ 
			and $m(x,\cdot)$ are locally Lipschitz continuous w.r.t. $u$ such that 
			\begin{align*}
				|l(t,x,u_1)-l(t,x,u_2)| & \leq C(1+|u_1|+|u_2|)|u_1-u_2| , \\ 
				|m(x,u_1)-m(x,u_2)| & \leq C(1+|u_1|+|u_2|)|u_1-u_2| \quad 
				\forall u_1,u_2\in \mathbb{R} , 
			\end{align*}
			for some constant $C>0$. 
			\item The initial condition $u_0:\Omega\rightarrow H$ is 
			$\mathcal{F}_0$- measurable with 
			\begin{align*}
				\E\left[\|u_0\|_H^p \right] <\infty,
			\end{align*}
			for all $p\geq 2$.
		\end{enumerate}
	\end{assumption}
	
	Under Assumption \ref{A1}, for any $\mathfrak{g}\in \mathbb{A}$, the equation \eqref{spde1} has a unique probabilistic strong solution in the variational setting, for the Gelfand triple
	\begin{align*}
		V\hookrightarrow H \hookrightarrow V',
	\end{align*}
	where $V=H^1$ is equipped with the norm $\|\cdot\|_{H^1}^2:=\|\cdot\|_{H}^2+\|\cdot\|_1^2$, see \cite{LR15} for further details.
	
	\begin{remark}
		For the stochastic Nagumo equation the non-linearity 
		\begin{align*}
			f(u)=\gamma u(u-1)(a-u)
		\end{align*}
		will only satisfy a one-sided Lipschitz condition, i.e. there exists $C>0$, such that for any $u,v\in \mathbb{R}$
		\begin{align*}
			\langle f(u)-f(v),u-v\rangle \leq C|u-v|^2,
		\end{align*}
		along with 
		\begin{align*}
			|f(u)|&\leq C(1+|u|^k),\\
			\langle f(u),u\rangle&\leq C(1+|u|^2),
		\end{align*}
		for any $u\in \mathbb{R}$ and some $k\in \mathbb{N}$. However our analysis can be adapted to this situation, which will be briefly discussed in Remark \ref{Rem1}, Remark \ref{Rem2}, Remark \ref{Rem3}, Remark \ref{Rem4} and Remark \ref{Rem5}.
	\end{remark}
	
	\subsection{Reduction to Feedback Controls}\label{reduction}
	
	In the whole paper we assume that for the control problem \ref{SCP} there exists an optimal control $\hat{\mathfrak{g}}\in \mathbb{A}$ of feedback type, i.e. 
	\begin{align*}
		\inf_{\mathfrak{g}\in \mathbb{A}}J(\mathfrak{g})=J(\hat{\mathfrak{g}})
	\end{align*}
	and
	\begin{align*}
		\hat{\mathfrak{g}}_t=\hat{G}(t,u_t^{\hat{G}}), 
	\end{align*}
	for some feedback map $\hat{G}:[0,T]\times H\rightarrow H$, such that the equation 
	\begin{equation}\label{spdeFeedback}
		\begin{cases}
			du^{\hat{G}}_t=[ A u^{\hat{G}}_t+\mathcal{F}(u^{\hat{G}}_t) 
			+\hat{G}(t,u_t^{\hat{G}})]dt+ B dW_t, \quad t\in [0,T]\\
			u^{\hat{G}}_0=u_0,
		\end{cases}
	\end{equation} 
	has a unique probabilistic strong solution, in the variational setting. Furthermore we will assume that the feedback map $\hat{G}$ is jointly continuous and satisfies a linear growth condition
	\begin{align*}
		\| \hat{G}(t,u) \|_H & \leq C(1+\| u \|_H).
	\end{align*}
	
	The simplest example where our assumptions are satisfied is the case of a linear quadratic control problem, where $\mathcal{F}$ is linear and $l(t,x,u)=|u|^2$, $m(t,u)=|u|^2$. In this case one can derive an explicit representation for the optimal control in terms of the solution to a Riccati equation, see e.g. \cite{Tud90}. In the case of the controlled stochastic heat equation, i.e. $\mathcal{F}\equiv 0$, the optimal control $\mathfrak{g}^{\ast}$ is given by
	\begin{equation}
		\mathfrak{g}^{\ast}_t = P(t) u^{\mathfrak{g}^{\ast}}_t,
	\end{equation}
	where $P:[0,T]\to L(H)$ is the solution of the Riccati equation (see e.g. \cite{Tud90})
	\begin{equation}\label{Ricc}
		\begin{cases}
			\partial_t P(t) + P(t) \Delta +\Delta P(t) - \text{Id} + P^2(t) = 0,\;\; t\in [0,T]\\
			P(T) = - \text{Id}.
		\end{cases}
	\end{equation}
	We will use this example later as a benchmark for our algorithm.
	
	In the following we will give a short explanation how one can construct such an optimal 
	feedback function $\hat{G}$ in a more general situation. To this end we consider the 
	Hamilton-Jacobi-Bellman equation 
	\begin{equation}\label{HJB}
		\begin{cases}
			\partial_t V + \frac{1}{2}\text{tr}(BQB^\ast D^2V)  
			+ \langle DV,Au+\mathcal{F}(u)\rangle_H 
			+ \int_\Lambda l(t,x,u(x)) \, dx \\
			\quad\qquad + \inf_{G\in \mathcal{U}} 
			\{\langle Dv,G\rangle + \|G\|_H^2\} =0, 
			\quad u\in H, t\in [0,T]   \\ 
			V_T(u)=\int_\Lambda m(x,u(x))\, dx , u \in H.
		\end{cases}
	\end{equation} 
	We assume that the HJB equation (\ref{HJB}) has a unique mild solution 
	$V:[0,T]\times H\rightarrow \mathbb{R}$ in the sense of 
	\cite[Definition~4.70]{FGS17a}, such that
	\begin{align*}
		V & \in \mathcal{C}^{0,1}([0,T]\times H,\mathbb{R}).
	\end{align*}
	If in addition the function $\gamma(p) := \text{arginf}_{g\in \mathcal{U}} 
	\{\langle p,g\rangle+\|g\|_{H}^2\}$ is continuous and equation 
	(\ref{spdeFeedback}) has a unique strong solution for  
	$\hat{G}(t,u) = \gamma( DV(t,u))$, then an optimal control of the control Problem \ref{SCP} is given by
	\begin{align*}
		\hat{\mathfrak{g}}_t & = \hat{G}(t,u_t^{\hat{G}}).
	\end{align*}
	For the existence, regularity of a solution to the HJB equation (\ref{HJB}) 
	and examples we refer to \cite{FGS17a}, in particular 
	Section 4.8.1 and Section 6.11. A sufficient condition for \eqref{spdeFeedback} to have a unique strong solution is that $\hat{G}$ is Lipschitz continuous in $u$, see e.g. \cite{LR15}. This is in particular the case if the solution $V$ of the HJB equation \eqref{HJB} has a bounded second derivative in $u$, see \cite[Theorem~4.155,Remark 4.202]{FGS17a}. Note that $\hat{G}(t,u)$ is in particular 
	bounded if $\mathcal{U}$ is a bounded subset of $H$. 
	
	In view of the previous discussion, we consider the following feedback control problem.
	
	\begin{problem}[FCP]\label{FCP}
		Minimize 
		\begin{align*}
			J(G) := J_1(u^G) + J_2(G(\cdot , u_\cdot^G)),
		\end{align*}
		subject to the spde
		\begin{equation}\label{spdeFeedback2}
			\begin{cases}
				du^{G}_t & =[ A u^{G}_t  + \mathcal{F}(u^{G}_t) 
				+ G(t,u_t^{G})]dt + B dW_t, \quad t\in [0,T]\\
				u^{G}_0 & =u_0,
			\end{cases}
		\end{equation}
		over the set of admissible feedback controls
		\begin{align*}
			U_{\mathrm{ad}} := \{G\in \mathcal{C}([0,T]\times H,\mathcal{U})|\text{ eq. }(\ref{spdeFeedback2}) \text{ has a unique strong solution and } G(\cdot,u_{\cdot}^G) \in\mathbb{A}\}. 
		\end{align*} 
	\end{problem}
	
	\begin{remark}
		If $G\in \mathcal{C}([0,T]\times H,\mathcal{U})$ is Lipschitz continuous in the second variable, uniformly in $t$, then equation \eqref{spdeFeedback} has a unique probabilistic strong solution, see again \cite{LR15} for details. In particular it holds $G\in U_{ad}$. Furthermore it is obvious that under the above assumption we have 
		\begin{align*}
			\inf_{\mathfrak{g}\in\mathbb{A}}J(\mathfrak{g})=\inf_{G\in U_{\mathrm{ad}}}J(G).
		\end{align*}
		With an abuse of notation, we use the same notation $J$ for the cost functional of the control problem \ref{SCP} and the feedback control problem \ref{FCP}.
	\end{remark}
	
	\subsection{Approximation Results}
	
	For the approximation of a feedback 
	function $G$, we consider a family of finite 
	dimensional subspaces $(S_h)_{h\in (0,1]}$ of $H^1$ satisfying 
	\begin{align*}
		\|P_hu-u\|_H\rightarrow 0, h\to 0, 
	\end{align*}
	where $P_h : H\rightarrow S_h$ denotes the $L^2$-orthogonal projection 
	onto $S_h$. We will also need the orthogonal projection 
	$\mathcal{R}_h:H^1\rightarrow S_h$ onto $S_h$ w.r.t. the inner product 
	\begin{align*}
		a(u,v) := \langle (-A)^{\frac{1}{2}} u,(-A)^{\frac{1}{2}} v\rangle_{H}.
	\end{align*}
	For any map $G:[0,T]\times H\rightarrow H$ we define the corresponding finitely based approximation of order $h\in (0,1]$ by 
	\begin{align}\label{FEApproximation}
		G^h:[0,T]\times H\rightarrow S_h, 
		G^h(t,u) := P_hG(t,P_hu).
	\end{align}
	Here finitely-based means that $G^h(t,u)$ depends on finitely many 
	coordinates of $u$ only.
	
	\subsubsection{Uniform Approximation Result}
	
	Our first main result provides sufficient conditions for a set of controls $\mathbb{U}\subseteq U_{ad}$ to lead to the optimal cost. This provides the theoretical foundation of our method.
	
	\begin{definition}\label{def:universal}
		Let $G:[0,T]\times H\rightarrow H$. We say that a subset $\mathbb{U}\subset U_{\text{ad}}$ satisfies the uniform approximation property with respect to $G$, if there exists a family $(G^{h,m})_{h\in (0,1],m\in \mathbb{N}}\subset \mathbb{U}$ that satisfies a linear growth condition uniformly in $m$, i.e. 
		\begin{align}\label{u1}
			\|G^{h,m}(t,u)\|_{H}&\leq C_h(1+\|u\|_{H}),
		\end{align}
		for some constant $C_h>0$, such that for any $R>0$
		\begin{align}\label{u2}
			\lim\limits_{m\rightarrow \infty}\sup_{(t,u)\in [0,T]\times \mathcal{B}_{H}(0,R)}\|G^h(t,u)-G^{h,m}(t,u)\|_{H}^2 = 0,
		\end{align}
		where 
		\begin{align}
			\mathcal{B}_{H}(0,R):=\{u\in H|\|u\|_{H}&\leq R\}
		\end{align}
		and $G^h$ is given as in \eqref{FEApproximation}.
	\end{definition}

	\begin{theorem}\label{main1}
		Assume that Assumption \ref{A1} is in force. Let $\mathbb{U}\subseteq U_{\mathrm{ad}}$ satisfy the uniform approximation property with respect to $\hat{G}$, then 
		\begin{align*}
			\inf_{\mathfrak{g}\in \mathbb{A}} J(\mathfrak{g}) 
			= \inf_{G\in \mathbb{U}}J(G).
		\end{align*}
	\end{theorem}

	\begin{remark}
		Theorem \ref{main1} tells us, if any finitely based approximation of 
		$\hat{G}$ can be approximated on bounded sets by elements of 
		$\mathbb{U}$, then we reach the optimal cost by minimizing over 
		$\mathbb{U}$. In Section \ref{ansatz} we will construct ansatz spaces $\mathbb{U}$ for our approximation method, that satisfy the uniform approximation property.
	\end{remark}

	\subsubsection{Convergence Rates}
	
	Sets of controls $\mathbb{U}\subseteq U_{ad}$ that satisfy the uniform approximation property with respect to the optimal control $\hat{G}$ are required to approximate every finitely based approximation $\hat{G}^h, h\in (0,1]$ on bounded sets. Therefore they might not be suited for explicit numerical implementation, since they typically involve controls with arbitrary high image dimension, as $h\downarrow 0$. In our second main theorem we are interested in error estimates for sets of controls that only approximate finitely based approximation $\hat{G}^h$ of a certain order $h\in (0,1]$. In order to formulate our second result we need to impose stronger assumptions on our control problem.
	
	\begin{assumption}\label{S1}
		\begin{enumerate}[S1.]
			\item The dispersion operator $B\sqrt{Q}$ takes values in $H^1$ 
			and is Hilbert-Schmidt, i.e. $\|B\sqrt{Q}\|_{L^0_{2,1}}<\infty$.
			\item The initial condition $u_0:\Omega\rightarrow H^1$ is 
			$\mathcal{F}_0$- measurable with 
			\begin{align*}
				\E\left[\|u_0\|_{H^1}^p\right] < \infty,
			\end{align*}
			for all $p\geq 2$.
			\item The admissible controls take values in $H^1$, in particular 
			we assume that
			\begin{align*}
				\mathcal{U} \subseteq H^1.
			\end{align*}
			Furthermore, the optimal feedback $\hat{G} (t,u):[0,T]\times 
			H\rightarrow \mathcal{U}$ is Lipschitz continuous in $u$ with 
			Lipschitz constant $L$ independent of $t$, i.e. 
			\begin{align*}
				\|\hat{G}(t,u_1)-\hat{G}(t,u_2)\|_H & \leq L \|u_1-u_2\|_H 
				\quad \forall t\in [0,T],u_1,u_2\in H 
			\end{align*}
		\end{enumerate}
	\end{assumption}
	
	Assumption \ref{S1} is still satisfied in the linear quadratic case (see Section \ref{reduction}). In a more general situation however, the existence of a solution to the HJB equation (\ref{HJB}) with Lipschitz regularity is difficult to check. We refer to Section \ref{ansatzFix} for a short discussion.
	
	In order to quantify the error resulting form the finitely based approximations $G^h$, we will specify the assumptions on our given finite element approximation: 
	
	\begin{assumption}\label{S2}
		\begin{enumerate}[R1.]
			\item For all $u\in 
			H^s,s\in \{1,2\}$ and $h\in (0,1]$ it holds
			\begin{align*}
				\|\mathcal{R}_hu-u\|_{H}\leq Ch^s\|u\|_{s}.
			\end{align*}
			Furthermore, for every $u\in H^1,h\in (0,1]$ it holds
			\begin{align*}
				\|P_h u\|_1\leq C\|u\|_{1}.
			\end{align*}
			\item The Ritz-projection coincides with the $L^2$-orthogonal projection on $H^1$, i.e. 
			\begin{align*}
				P_h|_{H^1}=\mathcal{R}_h
			\end{align*}
		\end{enumerate}
	\end{assumption}
	
	\begin{remark}
		\label{Projection}
		By the best approximation property of the orthogonal $L^2$-projection 
		we get for any $u\in H^s$, with $s\in \{1,2\}$
		\begin{align*}
			\|P_h u-u\|_{H} & \leq \|\mathcal{R}_h u - u\|_{H} 
			\leq C h^s\|u\|_{s}.
		\end{align*}
		Assumption R2 is satisfied for example when $A=\Delta$, where $\Delta$ 
		denotes the Laplace operator with Dirichlet boundary conditions on the unit 
		interval $(0,1)\subset \mathbb{R}$, and the finite dimensional subspaces 
		$S_h$ are given by
		\begin{align*} 
			S_h = \text{span}\{e_k,k=1,...,N\},
		\end{align*}
		for $h=\lambda_{N+1}^{-\frac{1}{2}}$, $N\in \mathbb{N}$, where 
		$e_k(\cdot)=\sqrt{2}\sin(k\pi\cdot)$ are the orthonormal eigenfunctions 
		of $A$ with corresponding eigenvalues $-\lambda_k=-k^2\pi^2$ (see 
		\cite{Kru14}).
	\end{remark}
	
	In our second main result we also want to cover the approximation error resulting from the spacial discretization of the controlled state equation. Therefore we introduce the following approximating control problem.
	
	\begin{problem}[FESCP]\label{FESCP}
		Let $h\in (0,1]$. Minimize 
		\begin{align*}
			J^h(G) = J_1(u^{G,h})+J_2(G (\cdot , u^{G,h})),
		\end{align*}
		over the set
		\begin{align*}
			U_{\mathrm{ad}}^L := \{G\in \mathcal{C}([0,T]\times H, 
			\mathcal{U}) |  G(t,u) \text{ is Lipschitz continuous in $u$, 
				uniformly in }t\},
		\end{align*}
		subject to the discretized SPDE
		\begin{equation}
			\label{FEspde}
			\begin{cases}
				\mathrm{d}u^{G,h}_t & = [A_h u^{G,h}_t  + P_h\mathcal{F}(u^{G,h}_t) 
				+ P_h G(t,u^{G,h}_t)] dt  + P_h B dW_t, \quad t\in [0,T] \\
				u^{G,h}_0 & = P_h u\in S_h,
			\end{cases}
		\end{equation}
		where $A_h u_h$ for $u_h\in S_h$ is defined as the unique element in 
		$S_h$ with 
		\begin{align*}
			a(u_h,v_h) = -\langle A_h u_h, v_h \rangle_{H}, \text{ for all }v_h\in S_h.
		\end{align*}
	\end{problem}
	
	The following lemma is easily shown by standard arguments.
	
	\begin{lemma} 
		\label{uniformest}
		For any $G\in U_{\mathrm{ad}}^L$ there exists a unique strong solution 
		$u^{G,h}:[0,T]\times \Omega\rightarrow S_h$ to \eqref{FEspde} and for any $p\ge 2$ it holds
		\begin{align*}
			\sup_{h\in (0,1]}\left(\sup_{t\in [0,T]} \mathbb{E} \left[ 
			\|u^{G,h}_t\|_H^p + \left( \int_0^T \|u^{G,h}_t \|_{H^1}^2 dt 
			\right)^{p/2} \right] \right) < \infty.
		\end{align*}
	\end{lemma}
	
	\begin{definition}\label{unifLipProp}
		Let $G:[0,T]\times H\rightarrow H$ and $h\in (0,1]$. We say that a sequence of subsets $(\mathbb{U}^{h,m})_{m\in \mathbb{N}}\subset U_{\text{ad}}^L$ satisfies the uniform Lipschitz approximation property of order $h$ with respect to $G$, if there exists a sequence of Lipschitz continuous controls $(G^{h,m})_{m\in \mathbb{N}}$ with Lipschitz constants independent of $m$, such that $G^{h,m}\in \mathbb{U}^{h,m}$, and for all $R>0$
		\begin{equation}
			\begin{aligned}
				\epsilon_m^{h,R} &   :=\sup_{(t,u)\in [0,T]\times \mathcal{B}_{H}(0,R)}\|P_h(G^h(t,u)-G^{h,m}(t,u))\|_{H}^2 	\label{RateEpsilon}  \rightarrow 0,
			\end{aligned}
		\end{equation}
		as $m\rightarrow \infty$.
	\end{definition}
	
	Our second main result will be split into two parts:
	
	\begin{theorem}\label{main2}
		Let $h\in (0,1]$. We assume that Assumption \ref{A1}, Assumption \ref{S1} and Assumption \ref{S2} R1 are in force. Let $(\mathbb{U}^{h,m})_{m\in \mathbb{N}}\subset U_{\text{ad}}^L$ satisfy the uniform Lipschitz approximation property of order $h$ with respect to $\hat{G}$, then it holds 
		\begin{align*}
			&\inf_{G\in \mathbb{U}^{h,m}}J^h(G)-\inf_{\mathfrak{g} 
				\in \mathbb{A}}J(\mathfrak{g})\\
			&\leq C \left( 1+ \E\left[ \int_{0}^{T} 
			\|\hat{G}(t,u_t^{\hat{G}})\|_{H^1}^2 \, dt \right]^{1/2}  
			\right) h + C_h\left( \epsilon^{h,R}_m+\frac{1}{R} \right)^{1/2},
		\end{align*}
		for all $R>0$, some uniform constant $C>0$ and some constant $C_h>0$ that is independent of $m$.
	\end{theorem}
	
	Under additional convexity assumptions on the coefficients we can prove 
	a lower bound on the approximating costs. To this end consider the 
	following modified assumptions  
	\begin{assumption}\label{convex}
		\begin{enumerate}[H1'.]
			\item The function $f:\mathbb{R}\rightarrow \mathbb{R}$ is linear.
			\item For any $t\in [0,T],x\in \Lambda$ and all $u_1,u_2 
			\in\mathbb{R}$ it holds
			\begin{align*}
				l(t,x,u_2)-l(t,x,u_1)-\partial_u l(t,x,u_1)(u_2-u_1)&\geq 0\\
				m(x,u_2)-m(x,u_1)-\partial_u m(x,u_1)(u_2-u_1)&\geq 0.
			\end{align*}
		\end{enumerate}
	\end{assumption}

	\begin{theorem}\label{main3}
		In addition to the assumptions of Theorem \ref{main2} we assume that Assumption \ref{S2} and Assumption \ref{convex} are in force. Then for any $\mathbb{U}\subseteq U_{\mathrm{ad}}$ it holds 
		\begin{align*}
			\inf_{\mathfrak{g}\in \mathbb{A}}J(\mathfrak{g})-\inf_{G\in \mathbb{U}}J^h(G)\leq Ch,
		\end{align*}
		for some constant $C>0$, which is independent of $h$.
	\end{theorem}
	
	\begin{remark}
		If $J^h$ is replaced by $J$ in Theorem \ref{main2} and Theorem \ref{main3} respectively, the obtained error bound still remains the same. Our proof can be easily adapted to this situation. 
	\end{remark}
	
	\section{Proof of Theorem \ref{main1}}\label{thm1}
	
	To simplify notations we denote $V=H^1$ equipped with the norm
	\begin{align*}
		\|\cdot\|_V:=\|\cdot\|_{H^1},
	\end{align*}
	$P_n := P_{\frac{1}{n}}$,  
	$G^n := \hat{G}^{\frac{1}{n}}$ and $G^{n,m}(t,u) :=  
	G^{\frac{1}{n},m}(t,u)$. Furthermore let 
	\begin{align*}
		u^n & :=u^{G^n}
	\end{align*}
	be the unique strong solution to equation (\ref{spdeFeedback2}) w.r.t. 
	$G^n$ and 
	\begin{align*}
		u^{n,m} := u^{G^{n,m}}
	\end{align*}
	be the unique strong solution to equation (\ref{spdeFeedback2}) w.r.t. 
	$G^{n,m}$. Due to the linear growth assumptions on $\hat{G}$ and $G^{h,m}$, 
	one can obtain the following a-priori estimates by standard arguments.
	
	\begin{lemma} 
		\label{APrioriBound1}
		For any $p\geq 2$ it holds
		\begin{align*}
			\sup_{n\in \mathbb{N}} \mathbb{E} \left[\sup_{t\in [0,T]} 
			\|u^n_t\|_{H}^p + \left( \int_{0}^{T} \|u^n_t\|_{V}^2 dt\right)^{p/2}
			\right] 
			& \leq C(1+\|u_0\|_{H}^p).
		\end{align*}
		Furthermore for any $p\geq 2$ and fixed $n\in \mathbb{N}$ we have 
		\begin{align*}
			\sup_{m\in \mathbb{N}} \mathbb{E} \left[\sup_{t\in [0,T]}  
			\|u^{n,m}_t\|_{H}^p + \left( \int_{0}^{T} \|u^{n,m}_t\|_{V}^2 dt 
			\right)^{p/2} \right] 
			& \leq C_n(1+\|u_0\|_{H}^p).
		\end{align*}
	\end{lemma}
	
	\begin{remark}\label{Rem1}
		The a-priori estimate in Lemma \ref{APrioriBound1} can also be obtained, if $f$ only satisfies a monotone growth condition, since we only need the Nemytskii operator $\mathcal{F}$ to satisfy
		\begin{align*}
			\langle \mathcal{F}(u),u\rangle_H\leq C(1+\|u\|_H^2),
		\end{align*}
		for every $u\in H$. 
	\end{remark}
	
	The following tightness result is a standard consequence of the a-priori bound given in Lemma \ref{APrioriBound1}, for further details see \cite{FG95} or \cite{RSZ22}. 
	
	\begin{lemma}
		\label{TightnessLemma}
		The sequence $(u^n)_{n\in \mathbb{N}}$ is tight in $\mathcal{C} ([0,T],V') 
		\cap L^2([0,T],H)$ and w.r.t. the weak topology in $L^2([0,T],V)$. Furthermore for any fixed $n\in \mathbb{N}$ the sequence $(u^{n,m})_{m\in \mathbb{N}}$ 
		is tight in $\mathcal{C}([0,T],V')\cap L^2([0,T],H)$ and w.r.t. the weak 
		topology in $L^2([0,T],V)$.
	\end{lemma}
	
	\subsection{Finitely Based Approximation of the Optimal Cost} 
	\label{finitelyBasedproof}
	
	In the first part we will prove that the optimal cost can be approximated 
	by the cost of the finitely based vector-fields 
	$\mathfrak{g}^n_t := G^{n}(t,u^n_t)$. In particular we will show that for 
	\begin{equation}
		\label{main1R1}
		J(\mathfrak{g}^n) \rightarrow J(\hat{\mathfrak{g}}) 
		= \inf_{\mathfrak{g} \in \mathbb{A}}J(\mathfrak{g})
	\end{equation}
	along a subsequence. In the following we set 
	\begin{align*}
		\Gamma := (L^2([0,T],H)\cap \mathcal{C}([0,T],V'))\times \mathcal{C}
		([0,T],U_1),
	\end{align*}
	where $U_1$ is a Hilbert space such that the embedding $U\subset U_1$ is 
	Hilbert-Schmidt. 
	
	By Lemma \ref{TightnessLemma}, Prohorov's theorem and the Skorohod representation theorem, there exist a probability space $(\tilde{\Omega},\tilde{\mathcal{F}}, \tilde{\mathbb{P}})$ and $\Gamma$-valued 
	random variables $(\tilde{u}^n,\tilde{W}^n)_{n\in \mathbb{N}}$ and $(\tilde{u},\tilde{W})$, such that 
	\begin{enumerate}
		\item $\tilde{W}^n=\tilde{W}$ for any $n\in \mathbb{N}$, 
		$\tilde{\mathbb{P}}$-a.s.
		\item $\mathcal{L}(\tilde{u}^n,\tilde{W}^n)=\mathcal{L}(u^n,W)$
		\item it holds $\tilde{\mathbb{P}}$-a.s.,
		\begin{align*}
			\|\tilde{u}^n-\tilde{u}\|_{L^2([0,T],H)}+\|\tilde{u}^n-
			\tilde{u}\|_{\mathcal{C}([0,T],V')}\rightarrow 0.
		\end{align*}
	\end{enumerate}
	
	Thanks to the a-priori bound in Lemma \ref{APrioriBound1}, it is not difficult to see that $(\tilde{u}^n)_{n\in \mathbb{N}}$ satisfies 
	for any $p\geq 2$
	\begin{align*}
		\sup_{n\in \mathbb{N}} \tilde{\mathbb{E}} \left[ \sup_{t\in [0,T]} 
		\|\tilde{u}^n_t\|_H^p  
		+ \sup_{n\in \mathbb{N}} \left( \int_0^T  \| 
		\tilde{u}^n_t\|_V^2 dt  \right)^{p/2}\right] < \infty.
	\end{align*}
	Furthermore, since $\|\tilde{u}^n-\tilde{u}\|_{L^2 ([0,T],H)} 
	\rightarrow 0$, $\tilde{	\mathbb{P}}$-a.s. , we have by using Fatou's lemma 
	\begin{align*}
		& \tilde{\E} \left[ \sup_{t\in [0,T]} \|\tilde{u}_t\|_H^p  
		+ \left(  \int_0^T  \|\tilde{u}_t\|_V^2 dt \right)^{p/2} \right] \\
		& \leq \liminf_{n\rightarrow\infty} \tilde{\E} \left[ \sup_{t\in [0,T]} 
		\|\tilde{u}_t^n\|_H^p + \left( \int_0^T  \|\tilde{u}_t^n\|_V^2 dt 
		\right)^{p/2}\right]  < \infty.
	\end{align*} 
	This implies in particular, that $\tilde{u}\in L^2 \left( \tilde{\Omega} , 
	\tilde{\mathbb{P}}, L^2 ([0,T], V)\right)$ and $\tilde{u}^n \rightarrow 
	\tilde{u}$ weakly. 
	
	We will now prove that $\tilde{u}$ is a weak solution to 
	\eqref{spdeFeedback2}. Recall that 
	\begin{align*}
		\langle Au,v\rangle = -\langle (-A)^{\frac{1}{2}} u,(-A)^{\frac{1}{2}} v\rangle_H , 
	\end{align*}
	so that $\|Au\|_{V'}\le \|u\|_V$, which implies that 
	\begin{align*}
		\tilde{\E}\left[\int_0^T \|A \tilde{u}^n_t\|_{V'}^2 dt \right] < \infty.
	\end{align*}
	For any $v\in L^2 \left( \tilde{\Omega} , \tilde{\mathbb{P}}, L^2 ([0,T], V)\right)$ it follows that 
	\begin{align*}
		\tilde{\E}\left( \int_0^T \langle A\tilde{u}_t^n , v_t\rangle \, dt \right) 
		& = \tilde{\E}\left( \int_0^T \langle \tilde{u}_t^n , Av_t\rangle \, dt \right) \\
		& \rightarrow \tilde{\E}\left( \int_0^T \langle \tilde{u}_t , Av_t\rangle \, dt \right) 
		= \tilde{\E}\left( \int_0^T \langle A\tilde{u}_t , v_t\rangle \, dt \right) 
	\end{align*}
	which implies that $A\tilde{u}^n \rightarrow A\tilde{u}$ 
	weakly in $L^2 \left( \tilde{\Omega} , \tilde{\mathbb{P}}, L^2 ([0,T], V')\right)$. 
	
	It remains to investigate the nonlinear drift coefficients. 
	
	\begin{lemma}
		\label{Lem1}
		There exists a subsequence $(\tilde{u}^{n_k})_{k\in \mathbb{N}}$ such that  
		\begin{align*}
			\int_0^T \| G^{n_k} (t,\tilde{u}^{n_k}_t) - \hat{G}(t,\tilde{u}_t)\|_H^2 dt 
			\rightarrow 0 \quad \tilde{\mathbb{P}}\text{-a.s. and in } L^p (\tilde{\Omega}, \tilde{\mathbb{P}}),    
		\end{align*}
		for all $p\ge 2$ and 
		\begin{align*}
			\int_0^T \| \mathcal{F}(\tilde{u}^{n_k}_t)-\mathcal{F}(\tilde{u}_t)\|_H^2 dt \rightarrow 0 
			\quad \tilde{\mathbb{P}}\text{-a.s. and in } L^p (\tilde{\Omega}, \tilde{\mathbb{P}}),
		\end{align*}
		for all $p\ge 2$. Furthermore, the process $\tilde{u}$ is a $H$-valued continuous process that 
		satisfies 
		\begin{align*}
			\langle \tilde{u}_t,v\rangle 
			= \langle u_0 + \int_0^t A\tilde{u}_t dt + \int_0^t \mathcal{F}(\tilde{u}_t) dt 
			+ \int_0^t \hat{G}(t,\tilde{u}_t) dt  + B\tilde{W}_t , v\rangle ,\quad  t\in [0,T], 
		\end{align*}
		for all $v\in V$, $\tilde{\mathbb{P}}$-a.s. as a $V'$-valued process.
	\end{lemma}
	
	\begin{proof}
		We just need to show the first statement of the lemma, the second part follows by standard arguments, see e.g. \cite{LR15}.
		We first observe that
		\begin{equation}
			\label{Eqn1Lem1}
			\begin{aligned} 
				&\| G^n(t,\tilde{u}_t^n)-\hat{G}(t,\tilde{u}_t)\|_H \\
				& = \| P_n \hat{G}(t,P_n \tilde{u}_t^n) - \hat{G}(t,\tilde{u}_t)\|_H \\ 
				& \le \|P_n \hat{G}(t,P_n \tilde{u}_t^n) - P_n \hat{G}(t, P_n \tilde{u}_t)\|_H  \\
				& \quad   + \|P_n \hat{G}(t,P_n \tilde{u}_t) - P_n \hat{G}(t, \tilde{u}_t)\|_H + \|P_n \hat{G}(t,\tilde{u}_t) - \hat{G}(t,\tilde{u}_t)\|_H  \\ 
				& \le \| \hat{G}(t,P_n \tilde{u}_t^n) - \hat{G}(t, P_n \tilde{u}_t)\|_H 
				+  \| \hat{G}(t,P_n \tilde{u}_t) - \hat{G}(t, \tilde{u}_t)\|_H \\
				&  \quad  + \|P_n \hat{G}(t,\tilde{u}_t) - \hat{G}(t,\tilde{u}_t)\|_H . 
			\end{aligned}
		\end{equation} 
		Since $\tilde{u}^n \to \tilde{u}$ in $L^2 ([0,T],H)$ $\tilde{	\mathbb{P}}$- a.s. it follows, passing to some subsequence again 
		denoted with $(\tilde{u}^n)_{n\in \mathbb{N}}$, that $\tilde{u}_t^n \to \tilde{u}_t$ in $H$ $\tilde{\mathbb{P}}\otimes dt$-a.s. . Now \eqref{Eqn1Lem1}, $\|P_n v-v\|_H\to 0$ for all $v\in H$ 
		and continuity of $\hat{G}$ imply that 
		\begin{align*}
			\| G^n(t,\tilde{u}_t^n)-\hat{G}(t,\tilde{u}_t)\|_H 
			\rightarrow 0\quad \tilde{\mathbb{P}}\otimes dt-a.e. .
		\end{align*}
		For convergence of $\int_0^T \| \hat{G}^n (t,\tilde{u}^n_t) - \hat{G}(t,\tilde{u}_t)\|_H^2 dt$ 
		in $L^p (\tilde{\Omega}, \tilde{\mathbb{P}})$ for all $p\ge 2$, it now suffices to show  
		\begin{equation} 
			\label{Eqn2Lem1}
			\sup_n \tilde{\mathbb{E}} \left[ \int_0^T \| G^n (t,\tilde{u}^n_t) - \hat{G}(t,\tilde{u}_t)\|_H^p dt \right] 
			< \infty, 
		\end{equation} 
		for all $p\ge 2$, since the latter implies uniform integrability of 
		$\| G^n (t,\tilde{u}^n_t) - \hat{G}(t,\tilde{u}_t)\|_H^p $, and thus 
		\begin{align*}
			\int_0^T \| G^{n_k} (t,\tilde{u}^{n_k}_t) - \hat{G}(t,\tilde{u}_t)\|_H^2 dt 
			\rightarrow 0 \quad \tilde{\mathbb{P}}\text{-a.s. and in } L^p (\tilde{\Omega}, \tilde{\mathbb{P}}),    
		\end{align*}
		for all $p\ge 2$.  But this follows from the a-priori bound in Lemma \ref{APrioriBound1} and the linear growth assumption on $\hat{G}$ by 
		\begin{align*} 
			\tilde{\mathbb{E}} \left[  \int_0^T \| G^n (t,\tilde{u}^n_t) - \hat{G}(t,\tilde{u}_t)\|_H^p dt \right] 
			& \le 2^p\, \tilde{\mathbb{E}} \left[  \int_0^T \| P_n \hat{G} (t,\tilde{u}^n_t)\|_H^p  
			+ \|\hat{G}(t,\tilde{u}_t)\|_H^p dt \right] \\ 
			& \le 2^p C\, \tilde{\mathbb{E}} \left[  \int_0^T 2 + \|\tilde{u}^n_t\|_H^p  
			+ \|\tilde{u}_t\|_H^p  dt \right] \\ 
			& \le 2^{p+1} C\, T \left( 1 + \sup_{n\in \mathbb{N}}\tilde{\mathbb{E}} \left[ \sup_{t\in[0,T]} \|\tilde{u}^n_t\|_H^p \right]  \right) < \infty, 
		\end{align*} 
		for all $p\ge 2$.
		
		The convergence  
		\begin{align*}
			\int_0^T \|\mathcal{F}(\tilde{u}^{n_k}_t)-\mathcal{F}(\tilde{u}_t)\|_H^2 \, dt \rightarrow 0   
		\end{align*}
		$\tilde{\mathbb{P}}$-a.s. and in $L^p (\tilde{\Omega}, \tilde{\mathbb{P}})$ 
		for all $p\ge 2$ follows by similar calculations, due to the (Lipschitz) continuity of $f$.
	\end{proof}
	
	Since $\tilde{W}$ is a cylindrical $Q$-Wiener process with respect to the filtration $\tilde{\mathcal{F}}_t$ generated by $\{\tilde{u}_s,\tilde{W}_s|s\leq t\}$, it follows from Lemma \ref{Lem1} that $(\tilde{u},\tilde{W})$ is a weak solution of equation \eqref{spdeFeedback2} for the control $\hat{G}$. Thanks to the pathwise uniqueness of equation \eqref{spdeFeedback2} for the control $\hat{G}$ we obtain $\mathcal{L}(\tilde{u}|\tilde{\mathbb{P}})=\mathcal{L}(u|\mathbb{P})$.
	
	\begin{remark}\label{Rem2}
		For the statement of Lemma \ref{Lem1} it is sufficient if the function $f$ is only continuous and one-sided Lipschitz continuous (so that the a-priori bound \ref{APrioriBound1} hold).
	\end{remark}
	
	\begin{lemma}\label{pointwise}
		There exists a subsequence $(\tilde{u}^{n_k})_{k\in \mathbb{N}}$ such that for any $t\in [0,T]$
		\begin{align*}
			\|\tilde{u}_t^{n_k}-\tilde{u}_t\|_H\rightarrow 0,
		\end{align*}
		$\tilde{\mathbb{P}}$-a.s. and in $L^p(\tilde{\Omega},\tilde{\mathbb{P}})$.
	\end{lemma}
	
	\begin{proof}
		We have for any $t\in [0,T]$ by Itô's formula \cite[Theorem~4.2.5]{LR15}
		\begin{align*}
			\|\tilde{u}_t^n-\tilde{u}_t\|_H^2&=2\int_{0}^{t}\langle A(\tilde{u}_s^n-\tilde{u}_s), \tilde{u}_s^n-\tilde{u}_s\rangle ds+\int_{0}^{t}\langle \mathcal{F}(\tilde{u}_s^n)-\mathcal{F}(\tilde{u}_s),\tilde{u}_s^n-\tilde{u}_s\rangle_Hds\\
			&\quad + \int_{0}^{t}\langle G^n(s,\tilde{u}_s^n)-\hat{G}(s,\tilde{u}_s),\tilde{u}_s^n-\tilde{u}_s\rangle_H ds\\
			&\leq -2\int_{0}^{t}\|\tilde{u}_s^n-\tilde{u}_s\|_V^2ds+2\int_{0}^{t}\|\tilde{u}_s^n-\tilde{u}_s\|_H^2ds\\
			&\quad +\int_{0}^{t}\|\mathcal{F}(\tilde{u}_s^n)-\mathcal{F}(\tilde{u}_s)\|_H^2ds + \int_{0}^{t}\|G^n(s,\tilde{u}_s^n)-\hat{G}(s,\tilde{u}_s)\|_H^2 ds.
		\end{align*}
		Using the Lipschitz continuity of $f$ and Gronwall inequality, we end up with 
		\begin{align*}
			\|\tilde{u}_t^n-\tilde{u}_t\|_H^2&\leq C_T\int_{0}^{T}\|G^n(s,\tilde{u}_s^n)-\hat{G}(s,\tilde{u}_s)\|_H^2 ds.
		\end{align*}
		Now by Lemma \ref{Lem1} there exists a subsequence, again denoted by $(\tilde{u}_t^n)_n$, such that 
		\begin{align*}
			\int_{0}^{T}\|G^n(s,\tilde{u}_s^n)-\hat{G}(s,\tilde{u}_s)\|_H^2 ds\rightarrow 0,
		\end{align*}
		$\tilde{\mathbb{P}}$-a.s. and in $L^p(\tilde{\Omega},\tilde{\mathbb{P}})$. This finishes the proof.
	\end{proof}
	
	\begin{remark}\label{Rem3}
		In the above proof of Lemma \ref{pointwise} we only need the one-sided Lipschitz continuity of $f$ to see that
		\begin{align*}
			\langle u-v,\mathcal{F}(u)-\mathcal{F}(v)\rangle_H\leq C\|u-v\|^2_H,
		\end{align*}
		for any $u,v\in H$.
	\end{remark}
	
	Given a probability space $(\Omega^\mu,\mathcal{F}^\mu,\mathbb{P}^\mu)$, we define for any $u\in L^2([0,T]\times \Omega^{\mu},H)$, $\mathfrak{g}\in L^2([0,T]\times \Omega^{\mu},H)$
	\begin{align*}
		J_1^\mu(u)  
		& := \E^\mu\left[\int_{0}^{T}\int_{\Lambda} l (t,x,u_t(x)) 
		\mathrm{d}x dt  + \int_{\Lambda}m(x,u_T(x))dx \right], \\
		J_2^\mu(\mathfrak{g}) 
		& := \E^\mu\left[\int_{0}^{T} \| \mathfrak{g}_t\|^2_H dt  
		\right].
	\end{align*}
	
	To show (\ref{main1R1}) we need the following lemma.
	
	\begin{lemma}
		\label{LipCost}
		Let $\mu$ be a probability space and $u_1,u_2\in L^2([0,T]\times \Omega^{\mu},H)$,
		$\mathfrak{g}_1,\mathfrak{g}_2\in L^2([0,T]\times \Omega^{\mu},H)$, such that 
		\begin{align*}
			\|\mathfrak{g}_1\|_{L^2([0,T]\times \Omega^{\mu},H)}+\|\mathfrak{g}_2\|_{L^2([0,T]\times \Omega^{\mu},H)}\leq \kappa,
		\end{align*}
		and
		\begin{align*}
			\|u_1\|_{L^2([0,T]\times \Omega^{\mu},H)}+\|u_2\|_{L^2([0,T]\times \Omega^{\mu},H)}\leq \kappa.
		\end{align*}
		Then it holds 
		\begin{align*}
			|J_1^{\mu}(u_1)-J_1^{\mu}(u_2)|\leq C(\kappa)\left[\|u_1-u_2\|_{L^2([0,T]\times \Omega^{\mu},H)} 
			+ \|u_1(T,\cdot)-u_2(T,\cdot)\|_{L^2(\Omega^\mu,H)}\right]
		\end{align*}
		and 
		\begin{align*}
			|J_2^{\mu}(\mathfrak{g}_1)-J_2^{\mu}(\mathfrak{g}_2)|\leq C(\kappa)\|\mathfrak{g}_1-\mathfrak{g}_2\|_{L^2([0,T]\times \Omega^{\mu},H)},
		\end{align*}
		for some constant $C(\kappa ) > 0$.
	\end{lemma}

	\begin{proof}
		Using the assumptions on $l$ and $m$ we get
		\begin{align*}
			& |J_1(u_1)-J_1(u_2)| \\
			& \leq \E^\mu \left[ \int_0^T \int_\Lambda |l (t,x,u_1(t,x)) 
			- l (t,x,u_2(t,x))| dx dt \right] \\
			& \quad + \E^\mu\left[ \int_\Lambda |m(x,u_1(T,x))-m(x,u_2(T,x))| 
			dx \right] \\
			&\leq C\, \E^\mu\left[ \int_0^T \int_\Lambda \left( 1 + |u_1(t,x)| 
			+ |u_2(t,x)| \right) |u_1(t,x)-u_2(t,x)| dx dt \right] \\
			& \quad + C\, \E^\mu \left[ \int_\Lambda \left( 1 + |u_1(T,x)|   
			+ |u_2(T,x)| \right) |u_1(T,x)-u_2(T,x)| dx \right].
		\end{align*}
		On the other hand we have 
		\begin{align*}
			& |J_2 (\mathfrak{g}_1) - J_2 (\mathfrak{g}_2) | \\
			& \leq \E^\mu\left[ \int_0^T \left| \|\mathfrak{g}_1(t)\|_H^2 - \|\mathfrak{g}_2(t)\|_H^2 \right| \, dt \right] \\
			& \leq \E^\mu\left[ \int_0^T \int_{\Lambda} (1 + |\mathfrak{g}_1(t,x)| + |\mathfrak{g}_2(t,x)|) 
			|\mathfrak{g}_1(t,x)-\mathfrak{g}_2(t,x)| dx dt \right].
		\end{align*}
		An application of the Cauchy-Schwarz inequality now yields the result.
	\end{proof}
	
	\begin{lemma}
		There exists a subsequence $(\tilde{u}^{n_k})_{k\in \mathbb{N}}$, such that
		\begin{align*}
			J(G^{n_k}(\cdot,u_{\cdot}^{n_k})) \rightarrow\inf_{\mathfrak{g}\in \mathbb{A}}J(\mathfrak{g})
		\end{align*}
	\end{lemma}
	
	\begin{proof}
		By Lemma \ref{Lem1} and Lemma \ref{pointwise}, there exists a subsequence $(\tilde{u}^{n_k})_{k\in \mathbb{N}}$, such that  
		\begin{align*}
			\|\tilde{u}^{n_k}-\tilde{u}\|_{L^2([0,T],H)}+\|\tilde{u}^{n_k}_T-\tilde{u}_T\|_{H}+\int_{0}^{T}\|G^{n_k}(t,\tilde{u}_t^{n_k})-\hat{G}(t,\tilde{u}_t)\|_H^2dt\rightarrow 0,
		\end{align*}
		$\tilde{\mathbb{P}}$-almost surely and in $L^p(\tilde{\Omega},\tilde{\mathbb{P}})$. Furthermore it holds 
		$\mathcal{L}(\tilde{u}|\tilde{\mathbb{P}})=\mathcal{L}(u|\mathbb{P})$ and $\mathcal{L}(\tilde{u}^n|
		\tilde{\mathbb{P}})=\mathcal{L}(u^n|\mathbb{P})$. Therefore by the a-priori estimates on  
		$\tilde{u}^n,\tilde{u}$ and the previous lemma, we have that for $\tilde{\mu}:=(\tilde{\Omega}, 
		\tilde{\mathcal{F}},\tilde{\mathbb{P}})$ it holds
		\begin{align*}
			J(G^{n_k}(\cdot,u_{\cdot}^{n_k}))&=J_1^{\tilde{\mu}}(\tilde{u}^{n_k})+ J_2^{\tilde{\mu}}(G^{n_k}(\cdot,\tilde{u}_{\cdot}^{n_k})) \\
			&\rightarrow J_1^{\tilde{\mu}}(\tilde{u})+J_2^{\tilde{\mu}}(\hat{G}(\cdot,\tilde{u}_{\cdot}))\\
			&=J(\hat{G}(\cdot,u^{\hat{G}}_{\cdot}))\\
			&=\inf_{\mathfrak{g}\in \mathbb{A}}J(\mathfrak{g}).
		\end{align*}
	\end{proof} 
	
	\subsection{Approximation of the Finitely Based Minimizing Sequence}
	
	We turn to the the second part of the proof for Theorem \ref{main1}. We will show that for any $n\in \mathbb{N}$ we have 
	\begin{align*}
		J(\mathfrak{g}^{n,m})\rightarrow J(\mathfrak{g}^n),
	\end{align*}
	along a subsequence, where $\mathfrak{g}^{n,m}_t:=G^{n,m}(t,u^{n,m}_t)$. The proof is quite similar to the previous proof. The only difference is the argument for the convergence result for $(\tilde{u}^{n,m})_{m\in \mathbb{N}}$, due to the different assumption on the approximating sequence $G^{n,m}$, $m\ge 1$. 
	
	Let us fix $n\in \mathbb{N}$. First note that due to the tightness of the sequence $(u^{n,m})_{m\in \mathbb{N}}$ stated in Lemma \ref{TightnessLemma} 
	there exist a probability space $(\tilde{\Omega},\tilde{\mathcal{F}},\tilde{\mathbb{P}})$ and a sequence of 
	$\Gamma$-valued random variables $(\tilde{u}^{n,m},\tilde{W}^m)_{m\in \mathbb{N}}$ and a $\Gamma$-valued random variable$(\tilde{u}^n,\tilde{W})$, such that 
	\begin{enumerate}
		\item $\tilde{W}^m = \tilde{W}$, $\tilde{\mathbb{P}}$-a.s. for any $m\in \mathbb{N}$,
		\item $\mathcal{L}(\tilde{u}^{n,m},\tilde{W}^m) = \mathcal{L}(u^{n,m},W)$, for any $m\in \mathbb{N}$
		\item it holds $\tilde{\mathbb{P}}$-a.s.,
		\begin{align*}
			\|\tilde{u}^{n,m}-\tilde{u}^n\|_{L^2([0,T],H)} + \|\tilde{u}^{n,m} - \tilde{u}^n \|_{\mathcal{C}([0,T],V')} 
			\rightarrow 0.
		\end{align*}
	\end{enumerate}
	
	Again, for any $p\geq 2$, it holds that 
	\begin{align}
		\label{APrioriBound2}
		\sup_{m\in\mathbb{N}} \tilde{\mathbb{E}} \left[ \sup_{t\in [0,T]} \|\tilde{u}^{n,m}_t\|_H^p + \left( \int_0^T  \|\tilde{u}^{n,m}_t\|_V^2 dt \right)^{p/2}\right] 
		< \infty, 
	\end{align}
	and therefore also  
	\begin{align*}
		\tilde{\E} \left[ \sup_{t\in [0,T]} \|\tilde{u}^n_t\|_H^p +  \left( \int_0^T  \|\tilde{u}^n_t\|_V^2 dt  
		\right)^{p/2}\right] < \infty.
	\end{align*} 
	Again, we can conclude from this that $\tilde{u}^n\in L^2 \left( \tilde{\Omega} , \tilde{\mathbb{P}},  
	L^2 ([0,T], V)\right)$ and $\tilde{u}^{n,m} \rightarrow \tilde{u}^n$ weakly.

	To identify $\tilde{u}^n$ as a weak solution to \eqref{spdeFeedback2} we again have that 
	$A\tilde{u}^{n,m} \rightarrow A\tilde{u}^n$ weakly in $L^2 \left( \tilde{\Omega} , 
	\tilde{\mathbb{P}}, L^2 ([0,T], V')\right)$, so that it remains to investigate 
	the nonlinear drift coefficients.

	\begin{lemma}
		\label{Lem2}
		There exists a subsequence $(\tilde{u}^{n,m_k})_{k\in \mathbb{N}}$ such that 
		\begin{align*}
			\lim\limits_{m\rightarrow \infty}\tilde{\E} \left[ \int_0^T \| G^{n,m_k}(t,\tilde{u}^{n,m_k}_t)-G^n(t,
			\tilde{u}_t^n)\|_H^2 dt\right] = 0\quad \tilde{\mathbb{P}}\text{-a.s. and in } L^p (\tilde{\Omega}, \tilde{\mathbb{P}})    ,
		\end{align*}
		for all $p\ge 2$ and  
		\begin{align*}
			\lim\limits_{m\rightarrow \infty}\tilde{\E}\left[\int_0^T \| \mathcal{F}(\tilde{u}^{n,m_K}_t) 
			-\mathcal{F}(\tilde{u}_t^n)\|_H^2 dt \right] = 0 \quad \tilde{\mathbb{P}}\text{-a.s. and in } L^p (\tilde{\Omega}, \tilde{\mathbb{P}}),
		\end{align*}
		for all $p\ge 2$. 
		Furthermore the process $\tilde{u}^n$ is a $H$-valued continuous process that satisfies 
		\begin{align*} 
			\langle \tilde{u}_t^n,v\rangle 
			= \langle u_0^n + \int_0^t A\tilde{u}^n_t dt + \int_0^t \mathcal{F}(\tilde{u}^n_t) dt 
			+ \int_0^t G^n(t,\tilde{u}^n_t) dt  + B\tilde{W}_t , v\rangle ,\quad t\in [0,T], 
		\end{align*}
		for all $v\in V$, $\tilde{\mathbb{P}}$-a.s. as a $V'$-valued process.
	\end{lemma}

	\begin{proof}
		Again, we only need to show the first statement of the theorem, the second part follows again by standard 
		arguments. This time we can estimate 
		\begin{equation}
			\label{Eqn1Lem9}
			\begin{aligned} 
				&\| G^{n,m}(t,\tilde{u}_t^{n,m}) - G^n (t,\tilde{u}^n_t)\|_H \\
				& \le \| G^{n,m} (t, \tilde{u}_t^{n,m}) - G^n (t,\tilde{u}^{n,m}_t)\|_H 
				+ \| G^n (t, \tilde{u}_t^{n,m}) - G^n (t,\tilde{u}^n_t)\|_H \\
				& = I_1^m(t) + I_2^m(t)  . 
			\end{aligned}
		\end{equation} 
		Passing to some subsequence (still denoted by $(\tilde{u}^{n,m})_{m\in \mathbb{N}})$, we have that  
		\begin{align*}
			\|P_n\tilde{u}^{n,m}_t-P_n\tilde{u}^n_t\|_H\leq \|\tilde{u}^{n,m}_t-\tilde{u}^n_t\|_H\rightarrow 0 \quad\tilde{\mathbb{P}}\otimes dt\text{-.a.e.}  
		\end{align*} 
		The continuity of $\hat{G}$ then implies that  
		\begin{align*}
			I_2^m(t) \leq \|\hat{G}(t ,P_n\tilde{u}^{n,m}_{t}) -\hat{G}(t,P_n\tilde{u}^n_{t})\|_{H} \rightarrow 0  
			\quad\tilde{\mathbb{P}}\otimes dt\text{-.a.e.}
		\end{align*}
		
		For the proof of convergence of $I_1^m$ to $0$ let 
		$R := \sup_{t\in [0,T]} \|\tilde{u}_t^{n}\|_H < \infty$ $\tilde{\mathbb{P}}$-a.s. 
		Then for $\tilde{\omega}\in\tilde{\Omega}$ and $t\in [0,T]$ with $R(\tilde{\omega}) < \infty$ and
		$\| \tilde{u}^{n,m}_t (\tilde{\omega}) -\tilde{u}^n_t (\tilde{\omega})\|_{H}\rightarrow 0$, it holds
		\begin{align*}
			\tilde{u}^{n,m}_t(\tilde{\omega})\in \mathcal{B}_H(0,\tilde{R}(\tilde{\omega})),
		\end{align*} 
		for all $m\in \mathbb{N}$, for some $\tilde{R}(\tilde{\omega})\leq 0$. Therefore the assumptions on $G^{n,m}$ now imply that  
		\begin{align*}
			\| G^{n,m }(t, \tilde{u}^{n,m}_t (\tilde{\omega} ))  
			- G^{n}(t , \tilde{u}^{n,m}_t (\tilde{\omega}))\|_{H} \rightarrow 0. 
		\end{align*}
		Since this is true for $\tilde{\mathbb{P}}\otimes dt$-a.e. $(\tilde{\omega},t)\in \tilde{\Omega}\times [0,T]$ we conclude that 
		\begin{align*}
			I_1^m(t) \rightarrow 0 \quad\tilde{\mathbb{P}}\otimes dt\text{-a.s.}
		\end{align*}
		The convergence of $\int_0^T \| G^{n,m} (t,\tilde{u}^{n,m}_t) - G^n (t,\tilde{u}^n_t)\|_H^2 dt$ 
		in $L^p (\tilde{\Omega}, \tilde{\mathbb{P}})$ for all $p\ge 2$, now follows similar to the corresponding 
		argument in Lemma \ref{Lem1} from the uniform integrability, implied by the bound 
		\begin{equation} 
			\label{Eqn2Lem2}
			\sup_m \tilde{\mathbb{E}} \left[ \int_0^T \| G^{n,m} (t,\tilde{u}^{n,m}_t) - G^n (t,\tilde{u}^n_t)\|_H^p dt \right] 
			< \infty 
		\end{equation} 
		for all $p\ge 2$. The bound \eqref{Eqn2Lem2} can be proven exactly in the same way as the proof of \eqref{Eqn2Lem1}, using 
		the uniform linear growth condition on $G^{n,m}$ w.r.t. $m$. 
		
		Finally, the convergence  
		\begin{align*}
			\int_0^T \|\mathcal{F}(\tilde{u}^{n,m}_t) - \mathcal{F}(\tilde{u}^n_t)\|_H^2 \, dt \rightarrow 0   
		\end{align*}
		$\tilde{\mathbb{P}}$-a.s. and in $L^p (\tilde{\Omega}, \tilde{\mathbb{P}})$ 
		for all $p\ge 2$ follows by similar calculations, due to the (Lipschitz) continuity of $f$. 
	\end{proof}

	Again, thanks to the pathwise uniqueness of equation \eqref{spdeFeedback2} for the control $G^n$ and Lemma \ref{Lem2}, we obtain $\mathcal{L}(\tilde{u}^n|\tilde{	\mathbb{P}})=\mathcal{L}(u^n|\mathbb{P})$. The following can be proven in a completely similar way as Lemma \ref{pointwise}.
	
	\begin{lemma}\label{pointwise2}
		There exists a subsequence $(\tilde{u}^{n,m_k})_{k\in \mathbb{N}}$ such that for any $t\in [0,T]$
		\begin{align*}
			\|\tilde{u}_t^{n,m_k}-\tilde{u}^n_t\|_H\rightarrow 0,
		\end{align*}
		$\tilde{\mathbb{P}}$-a.s. and in $L^p(\tilde{\Omega},\tilde{\mathbb{P}})$.
	\end{lemma}
	
	We are now in the position to finish the proof of our first main theorem of this chapter.
	
	\begin{proofThm}{\ref{main1}} By Lemma \ref{Lem2} and Lemma \ref{pointwise2}, there exists a subsequence $(\tilde{u}^{n,m_k})_{k\in \mathbb{N}}$, such that
		\begin{align*}
			\|\tilde{u}^{n,m_k}-\tilde{u}^n\|_{L^2([0,T],H)}+\|\tilde{u}^{n,m_k}_T-\tilde{u}^n_T\|_{H}+\int_{0}^{T}\|G^{n,m_k}(t,\tilde{u}_t^{n,m_k})-G^n(t,\tilde{u}_t^n)\|dt\rightarrow 0,
		\end{align*}
		$\tilde{\mathbb{P}}$-almost sure and in $L^p(\tilde{\Omega},\tilde{	\mathbb{P}})$. Furthermore it holds  $\mathcal{L}(\tilde{u}^n|\tilde{\mathbb{P}})=\mathcal{L}(u^n|\tilde{\mathbb{P}})$ and $\mathcal{L}(\tilde{u}^{n,m}|\tilde{\mathbb{P}})=\mathcal{L}(u^{n,m}|\tilde{\mathbb{P}})$. Therefore Lemma \ref{LipCost} implies for $\tilde{\mu}:=(\tilde{\Omega}, 
		\tilde{\mathcal{F}},(\tilde{\mathcal{F}}_t)_t,\tilde{\mathbb{P}})$
		\begin{align*}
			J(G^{n,m_k}(\cdot,u_{\cdot}^{n,m_k}))&=J_1^{\tilde{\mu}}(\Tilde{u}^{n,m_k})+J_2^{\Tilde{\mu}}(G^{n,m_k}(\cdot,\Tilde{u}_{\cdot}^{n,m_k}))\\
			&\rightarrow J_1^{\tilde{\mu}}(\Tilde{u}^n)+J_2^{\Tilde{\mu}}(G^n(\cdot,\Tilde{u}_{\cdot}^n))\\
			&=J(G^n(\cdot,u_{\cdot}^n)).
		\end{align*}
		By this observation, together with the result from Section \ref{finitelyBasedproof}, we can construct a sequence $(\overline{G}^n)_{n\in \mathbb{N}}\subseteq \mathbb{U}$, such that
		\begin{align*}
			J(\overline{G}^n(\cdot,u^{\overline{G}^n}_{\cdot}))\rightarrow J(\hat{G}(\cdot, u_{\cdot}^{\hat{G}}))=\inf_{\mathfrak{g}\in \mathbb{A}}J(\mathfrak{g}).
		\end{align*}
		Therefore we obtain
		\begin{align*}
			\inf_{G\in \mathbb{U}}J(G)&\leq J(\overline{G}^n)\\
			&=J(\overline{G}^n(\cdot,u^{\overline{G}^n}_{\cdot}))\rightarrow \inf_{\mathfrak{g}\in \mathbb{A}}J(\mathfrak{g}),
		\end{align*}
		hence 
		\begin{align*}
			\inf_{G\in \mathbb{U}}J(G)&\leq \inf_{\mathfrak{g}\in \mathbb{A}}J(\mathfrak{g}). 
		\end{align*}
		Since $\overline{G}^n(\cdot,u_{\cdot}^{\overline{G}^n})\in \mathbb{A}$ for all $n\in \mathbb{N}$, we obtain  
		\begin{align*}
			\inf_{G\in \mathbb{U}}J(G)&=\inf_{\mathfrak{g}\in \mathbb{A}}J(\mathfrak{g}). 
		\end{align*}
	\end{proofThm} 
	
	\section{Proof of Theorem \ref{main2}, \ref{main3}}\label{thm2}
	
	\subsection{Finite Element Approximation}
	
	We start with upper and lower bound estimates for the optimal cost of the finite element discretiation \ref{FESCP} of the control problem \ref{SCP}. 
	
	\subsubsection{Upper Bound}
	
	The main theorem concerning the upper bound of the approximating cost 
	functional $J^h$ is the following:
	
	\begin{theorem}
		\label{FEapprox}
		Let $h\in (0,1]$. Assume that Assumption \ref{A1} is in force, then it holds
		\begin{align*}
			\inf_{g\in \mathbb{A}}J(g)\geq J^h(\hat{G})-\epsilon_1(h),
		\end{align*}
		where 
		\begin{align*}
			\epsilon_1(h) = C(T) \left( 1 + \mathbb{E} \left[ 
			\int_0^T \|\hat{G}(t,u_t^{\hat{G}}) \|_{H^1}^2 \, dt\right]^{1/2} \right)h,
		\end{align*}
		for some $C(T)>0$, independent of $h$. This shows in particular that 
		\begin{align*}
			\inf_{\mathfrak{g}\in \mathbb{A}}J(\mathfrak{g}) 
			\geq \inf_{G\in U_{\mathrm{ad}}^L}J^h(G) - \epsilon_1(h),
		\end{align*}
		since $\hat{G}$ is an admissible feedback control.
	\end{theorem}
	
	We will now use the rest of this subsubsection to prove Theorem 
	\ref{FEapprox}. First we quantify the difference between the cost of a 
	control and its finite element approximation in terms of the parameter $h$.
	
	\begin{lemma} 
		\label{FECost}
		Let $G\in U_{\mathrm{ad}}^L$, $u^{G}$ be the unique solution to 
		(\ref{spdeFeedback2}) and $u^{G,h}$ be the unique solution to 
		(\ref{FEspde}), for $h\in (0,1]$. Then it holds 
		\begin{align*}
			|J_1(u^{G})-J_1(u^{{G},h})|\leq C(T)h
		\end{align*}
		and 
		\begin{align*}
			|J_2(G(\cdot,u_{\cdot}^{G}))-J_2(P_hG(\cdot,u_{\cdot}^{G,h}))| 
			\leq C(T) h \mathbb{E}\left[ \int_0^T \| G(t,u_t^{G})\|_{H^1}^2\, dt 
			\right]^{1/2},
		\end{align*}
		for some constant $C(T)>0$ which only depends on the Lipschitz constant 
		of $G$, but is independent of $h$.
	\end{lemma}
	
	\begin{proof}
		By Lemma \ref{uniformest} it holds 
		\begin{align*}
			\|u^{G,h}\|_{L^2([0,T]\times \Omega,H)}\leq C,
		\end{align*}
		for some constant $C > 0$ which may depend on the Lipschitz constant 
		of $G$, but is independent of $h$. Using $\|\P_h G(t,u_t^{G,h})\|_H 
		\le \|G(t,u_t^{G,h})\|_H$ and the fact that $G(t, u)$ is at most of 
		linear growth w.r.t. $u$, due to the Lipschitz property, we obtain
		\begin{align*}		
			\mathbb{E} \left[ \int_0^T \|P_h G(t,u_t^{G,h}) \|_H^2 \, dt \right]
			& \leq C \mathbb{E} \left[ \int_{0}^{T} (1 + \| u^{G,h}_t \|_H^2) dt \right] 
			\leq C,
		\end{align*}
		for some constant $C>0$ that depends only on the Lipschitz constant of 
		$G$. Now we can apply Lemma \ref{LipCost} to obtain
		\begin{align*}
			|J_1(u^{G})-J_1(u^{{G},h})|^2 
			& \leq C \E\left[ \int_0^T\|u^{{G}}_t-u^{{G},h}_t\|_H^2\right] 
			+ C\E\left[\|u^{{G}}_T-u^{{G},h}_T\|_H^2\right].
		\end{align*}
		Using Lemma \ref{LipCost} and the Lipschitz continuity of $G$ again, we get
		\begin{equation} 
			\label{FECost:eq}
			\begin{aligned}
				| J_2(G(\cdot, u_{\cdot}^{G})) - J_2 (P_h G(\cdot,u_{\cdot}^{G,h})) |^2 
				& \leq C\, \E\left[ \int_{0}^{T} \|{G}(t,u^{{G}}_t) - P_h{G}(t,u^{{G},h}_t) 
				\| _H^2 dt \right] \\
				& \leq C\, \E\left[ \int_{0}^{T} \|P_h({G}(t,u^{{G}}_t) - G(t,u^{G,h}_t)) 
				\|_H^2 dt \right] \\
				& \quad + C\, \E\left[ \int_{0}^{T} \|G(t,u^{G}_t)  
				- P_h G(t,u^{G}_t) \|_H^2 dt\right] \\
				& \leq C\, \E\left[\int_{0}^{T} \|u^{G}_t  
				- u^{G,h}_t \|_H^2 dt \right] \\
				& \quad + C\, \E\left[\int_{0}^{T} \|G(t,u^{G}_t) 
				- P_h G (t,u^{G}_t) \|_H^2 dt \right],
			\end{aligned}
		\end{equation} 
		where the constant $C$ may differ from line to line, but will be 
		independent of $h$. Since $G(t,u^{G}_t)\in H^1$, we get from Remark 
		\ref{Projection}
		\begin{align*}
			\|G(t,u^{G}_t) - P_h G (t,u^{G}_t) \|_H^2 & \leq h^2\|G(t,u_t^G)\|_{H^1}^2
		\end{align*}
		and from \cite[Corollary~7.2]{Kru14} we get
		\begin{align*}
			\E \left[ \int_{0}^{T} \|u^{G}_t -u^{G,h}_t \|_H^2dt \right] \leq C(T) h^2 ,
		\end{align*} 
		for some constant $C(T)$ independent of $h$. Inserting both estimates 
		into \eqref{FECost:eq} yields the desired result.
	\end{proof}
	
	\begin{remark}\label{Rem4}
		If $f$ is only one-sided Lipschitz continuous, we refer to \cite{CH19} for the corresponding estimate of the term 
		\begin{align*}
			\mathbb{E}\left[\int_{0}^{T}\|u^{G}_t -u^{G,h}_t \|_H^2dt\right]
		\end{align*}
		in the above proof. 
	\end{remark}
	
	\begin{proofThm}{\ref{FEapprox}}
		Theorem \ref{FEapprox} is now a simpel consequence of Lemma \ref{FECost} by
		\begin{align*}
			\inf_{\mathfrak{g}\in \mathbb{A}}J(\mathfrak{g}) 
			& = J_1(u^{{\hat{G}}}) + J_2(\hat{G}(\cdot,u_{\cdot}^{\hat{G}})) \\
			& \geq J_1(u^{\hat{G},h}) + J_2(P_h\hat{G}
			(\cdot,u_{\cdot}^{\hat{G},h})) - C(T) \left( 1  
			+ \mathbb{E} \left[ \int_{0}^{T} \| \hat{G} (t,u_t^{\hat{G}}) \|_{H^1}^2  
			\, dt \right] \right) h \\
			& = J^h ( \hat{G})  - C(T) \left( 1 + \mathbb{E} \left[ \int_{0}^{T}
			\| \hat{G}(t,u_t^{\hat{G}}) \|_{H^1}^2 dt \right] \right) h.
		\end{align*}
	\end{proofThm}
	
	\subsubsection{Lower Bound}
	
	Under the additional convexity assumptions specified in Assumption 
	\ref{convex} we can also prove the lower bound. The main result 
	concerning the lower bound of the approximating cost functional $J^h$ is 
	the following: 
	
	\begin{theorem}\label{convexEst}
		Let $h\in (0,1]$. In addition to the assumptions of Theorem  
		\ref{FEapprox} we assume that the Assumption \ref{convex} and Assumption \ref{S2} are in force. 
		Then for any $G\in U_{\mathrm{ad}}$ it holds
		\begin{align*}
			J^h (G) + Ch   \geq J(\hat{G}),
		\end{align*} 
		for some constant $C>0$ independent of $h$ and $G$. This shows in 
		particular that 
		\begin{align*}
			\inf_{\mathfrak{g}\in \mathbb{A}}J(\mathfrak{g}) 
			\leq \inf_{G\in U_{\mathrm{ad}}}J^h(G)+Ch.
		\end{align*}
	\end{theorem}
	
	The proof of the lower bound uses the necessary optimality condition 
	for the control problem \ref{SCP}. To this end let  
	$H(t, \cdot ) : H\times \mathcal{U} \times H 
	\rightarrow \mathbb R$
	\begin{align*} 
		H(t,u,\mathfrak{g} ,p) = \int_\Lambda l(t, x, u(x))\, dx  
		+ \|\mathfrak{g}\|^2_H + \langle p, \mathcal {F} (u) + \mathfrak{g} 
		\rangle_H  
	\end{align*} 
	denote the reduced Hamiltonian associated with the control problem \ref{SCP} and let 
	$\mathcal L (t,u) = \int_{\Lambda} l(t,x, u(x))dx$, $t\in [0,T]$, 
	$u\in H$ and respectively 
	$\mathcal{M}(u)=\int_{\Lambda} m(x,u(x))dx$. 
	By Assumption H1' and H2' the Hamiltonian is convex in the variables 
	$u,g$, i.e. 
	\begin{multline*} 
		H(t,u_2, \mathfrak{g}_2 , p ) 
		- H ( t,u_1, \mathfrak{g}_1,p) \\
		- \partial_u 
		H(t,u_1, \mathfrak{g}_1, p ) (u_2 - u_1) - \partial_g H(t,u_1, 
		\mathfrak{g}_1, p) (\mathfrak{g}_2 - \mathfrak{g}_1) \geq 0,
	\end{multline*}
	for any $t\in [0,T]$, $u_1,u_2\in H^1$, $\mathfrak{g}_1, \mathfrak{g}_2 
	\in\mathbb{A}$ and $p\in H$.
	
	\cite[Theorem~7.2]{SW21} now proves a stochastic maximum 
	principle for the case where the running cost $l(t,x,u)$ does not depend 
	on $t$. A straightforward generalization to the time dependent case now 
	yields the following
	
	\begin{theorem}
		There exist adapted processes $(q,p)$ with 
		\begin{align*}
			p\in L^2([0,T] \times \Omega,H^1) \cap L^2(\Omega, 
			\mathcal{C}([0,T],H))
		\end{align*}
		and 
		\begin{align*}
			q\in L^2([0,T] \times \Omega, L_2^0 )
		\end{align*}
		satisfying 
		\begin{equation}\label{adjoint1}
			\begin{cases}
				dp_t 
				& = - \left[ A p_t + D\mathcal{F} (u_t^{\hat{\mathfrak{g}}}) p_t 
				+ D\mathcal{L} (t,u^{\hat{\mathfrak{g}}}_t) \right] dt  
				+ q_t dW_t ,  \quad t\in [0,T] \\
				p_T
				& = D \mathcal{M}(u_T^{\hat{\mathfrak{g}}}),
			\end{cases}
		\end{equation} 
		such that
		\begin{align*}
			H(t, u^{\hat{\mathfrak{g}}}, v, p_t) 
			\geq H(t, u^{\hat{\mathfrak{g}}}, \hat{\mathfrak{g}}_t, p_t),
		\end{align*}
		for all $v\in\mathcal U$ and almost every $(t,\omega)\in [0,T]  
		\times \Omega$. In particular it holds 
		\begin{align*}
			\langle \partial_g H(t, u^{\hat{\mathfrak{g}}}, \hat{\mathfrak{g}}, 
			p_t),  v-\mathfrak{\hat{g}}_t \rangle_H \geq 0,
		\end{align*}
		for all $v\in\mathcal U$ and almost every $(t,\omega)\in [0,T] 
		\times \Omega$.
	\end{theorem}
	
	The proof of Theorem \ref{convexEst} is inspired by the technique of the proof of \cite[Theorem~6.16]{CD18b}, where the authors prove convergence of the optimal cost for finite player optimization problems towards the optimal cost of the limiting McKean-Vlasov control problem.
	
	\begin{proofThm}{\ref{convexEst}}
		We first rewrite 
		\begin{align*}
			J^h(G)-J(\hat{G}) & = J_1(u^{G,h}) - J_1(u^{\hat{G}}) 
			+ J_2 (P_h G(\cdot,u_{\cdot}^{{G,h}})) - J_2 ( \hat{G}
			(\cdot,u_{\cdot}^{\hat{G}})) \\
			& = \mathbb{E} \left[ \int_{0}^{T} \int_{\Lambda} 
			l(t,x,u_t^{G,h}(x)) - l(t,x,u_t^{\hat{G}}(x)) dx\, 
			dt \right] \\
			& \quad + \mathbb{E} \left[ \int_{\Lambda}   
			m(x,u_T^{G,h}(x)) - m(x,u_T^{\hat{G}}(x))  dx \right] \\
			& \quad + \mathbb{E} \left[ \int_{0}^{T}  \|P_h G(t,u_t^{G,h}) 
			\|_H^2 - \|\hat{G}(t,u_t^{\hat{G}})\|_H^2  dt \right].
		\end{align*}
		Let $(p,q)$ be the solution to the adjoint equation \eqref{adjoint1} w.r.t. 
		$\hat{\mathfrak{g}}_t = \hat{G} (t,u_t^{\hat{G}})$, then we can write 
		\begin{align*}
			J^h(G) - J(\hat{G}) & = T_1 + T_2,
		\end{align*}
		where
		\begin{align*} 
			T_1 
			& = \mathbb{E} \left[ \langle u_T^{G,h} -u_T^{\hat{G}}, p_T 
			\rangle_H \right] 
			+ \mathbb{E} \left[ \int_{0}^{T} 
			\int_{\Lambda}  l(t,x,u_t^{G,h}(x)) - l(t,x,u_t^{\hat{G}} 
			(x))  dx\, dt \right] \\
			& \quad + \mathbb{E} \left[ \int_{0}^{T} \| P_h G(t,u_t^{G,h}) 
			\|_H^2 - \|\hat{G}(t,u_t^{\hat{G}}) \|_H^2 dt \right] \\
			T_2 
			& = \mathbb{E} \left[ \int_{\Lambda} m(x,u_T^{G,h}(x)) -
			m( x, u_T^{\hat{G}} (x) ) dx \right] 
			- \mathbb{E} \left[ \langle u_T^{G,h} - u_t^{\hat{G}} , 
			D\mathcal{M} (u_T^{\hat{G}}) \rangle_H \right].
		\end{align*}
		We start by estimating the term $T_2$. By the convexity of $m$ in $u$, we get 
		\begin{align*}
			m(x,u')-m(x,u)-\partial_u m(x,u)(u'-u)\geq 0,
		\end{align*}
		for every $x\in\Lambda$, $u',u\in \mathbb{R}$, hence 
		\begin{align*}
			T_2 &\geq 0.
		\end{align*}
		For the other term $T_1$ we consider the equation for $u_t^{G,h}-u_t^{\hat{G}}$:
		\begin{equation}
			\begin{cases}
				d(u_t^{G,h}-u_t^{\hat{G}})=[A_h u^{{G},h}_t-Au_t^{\hat{G}}+(P_h\mathcal{F}(u^{{G},h}_t)-\mathcal{F}(u_t^{\hat{G}}))\\
				\quad\qquad\qquad\qquad\qquad \qquad \qquad \quad +(P_hG(t,u_t^{G,h})-\hat{G}(t,u_t^{\hat{G}}))]dt\\
				\qquad \qquad \qquad \quad + (P_h-I) B dW_t, \quad t\in [0,T]\\
				(u_0^{G,h}-u_0^{\hat{G}})=(P_h-I)u_0.
			\end{cases}
		\end{equation}
		By It\^o's formula \cite[Lemma~2.15]{Par21} we get 
		\begin{align*}
			d\langle u_t^{G,h}-u_t^{\hat{G}}, p_t\rangle_H 
			& = \langle u_t^{G,h}-u_t^{\hat{G}}, dp_t \rangle_H 
			+ \langle p_t,d(u_t^{G,h} - u_t^{\hat{G}}) 
			\rangle_H
			+ d\langle u^{G,h} - u^{\hat{G}}, p\rangle_t.
		\end{align*}
		Recall that $P_h|_{H^1} = \mathcal{R}_h$ by R2 (Assumption \ref{S2}). This 
		now implies that 
		\begin{align*}
			\langle A_h u^{G,h}_t, p_t \rangle_H 
			& = \langle A_h u^{G,h}_t, P_h p_t\rangle_H 
			= \langle A_h u^{G,h}_t, \mathcal{R}_h p_t\rangle_H \\
			& = -\langle (-A)^{\frac{1}{2}}u^{G,h}_t,(-A)^{\frac{1}{2}} \mathcal{R}_h 
			p_t \rangle_H 
			= -\langle (-A)^{\frac{1}{2}} u^{G,h}_t, (-A)^{\frac{1}{2}} p_t\rangle_H \\
			& = \langle Au^{G,h}_t,  p_t \rangle , 
		\end{align*}
		so that 
		\begin{align*}
			\langle A_h u_t^{G,h} - Au_t^{\hat{G}}  , p_t \rangle  
			& - \langle  Ap_t,u_t^{G,h}-u_t^{\hat{G}} \rangle  
			= 0.
		\end{align*}
		
		Now
		\begin{align*}
			d\langle u_t^{G,h}-u_t^{\hat{G}}, p_t \rangle_H  
			& = -\langle u_t^{G,h} - u_t^{\hat{G}}, D\mathcal{F}
			(u_t^{\hat{G}})p_t + D\mathcal{L} (t, u_t^{\hat{G}}) 
			\rangle_H dt \\
			& \quad + \langle p_t, (P_h\mathcal{F}(u^{{G},h}_t)-\mathcal{F}
			(u_t^{\hat{G}})) + (P_hG(t,u_t^{G,h}) - \hat{G} (t,u_t^{\hat{G}})) 
			\rangle_H dt \\
			& \quad + \langle q_t^*(u_t^{G,h} - u_t^{\hat{G}}), dW_t \rangle_H\\
			& \quad + \langle (P_h-I)^*p_t, BdW_t \rangle_H \\
			& \quad + \langle q_t \sqrt{Q},(P_h-I)B\sqrt{Q} 
			\rangle_{L_2^0} dt.
		\end{align*}
		Taking expectation and invoking the definition of the Hamiltonian, we 
		arrive at     
		\begin{align*}
			T_1 
			& = \mathbb{E}\left[\int_{0}^{T} H(t,u_t^{\hat{G}}, P_h 
			G(t,u_t^{G,h}), p_t) - H(t, u_t^{\hat{G}},  
			\hat{G}(t, u_t^{\hat{G}}), p_t) dt \right] \\ 
			& \quad + \mathbb{E} \left[ \int_{0}^{T} \langle  p_t, 
			P_h \mathcal{F} (u_t^{G,h}) - \mathcal{F} (u_t^{\hat{G}})  
			- D\mathcal{F} ( u_t^{\hat{G}} )(u_t^{G,h} - u_t^{\hat{G}})\rangle_H
			dt \right] \\ 
			& \quad + \mathbb{E} \left[ \int_{0}^{T} \int_\Lambda l(t,x,u_t^{G,h} 
			(x)) - l(t,x,u_t^{\hat{G}} (x)) - \partial_u l (t,x, u_t^{\hat{G}} (x))
			(u_t^{G,h} - u_t^{\hat{G}}) (x) dx\, dt \right] \\ 
			& \quad + \int_0^T \langle q_t \sqrt{Q},(P_h-I)B\sqrt{Q} 
			\rangle_{L_2^0} dt .
		\end{align*} 
		The maximum principle and convexity of $l$ in $u$ now imply that 
		\begin{align*}
			T_1 
			& \ge \mathbb{E} \left[ \int_{0}^{T} \langle  p_t, 
			P_h \mathcal{F} (u_t^{G,h}) - \mathcal{F} (u_t^{\hat{G}})  
			- D\mathcal{F} ( u_t^{\hat{G}} )(u_t^{G,h} - u_t^{\hat{G}})\rangle_H
			dt \right]\\
			&\quad  + \int_0^T \langle q_t \sqrt{Q},(P_h-I)B\sqrt{Q} \rangle_{L_2^0} dt .
		\end{align*} 
		Since $f$ is linear, hence $P_h\mathcal{F} (u_t^{G,h}) = \mathcal{F} ( P_h u_t^{G,h}) = \mathcal{F} (u_t^{G,h})$, the first term on the right hand side vanishes, i.e. 
		\begin{equation*}
			\mathbb{E} \left[ \int_{0}^{T} \langle  p_t, 
			P_h \mathcal{F} (u_t^{G,h}) - \mathcal{F} (u_t^{\hat{G}})  
			- D\mathcal{F} ( u_t^{\hat{G}} )(u_t^{G,h} - u_t^{\hat{G}})\rangle_H 
			dt \right] = 0.
		\end{equation*}
		Furthermore Assumption \ref{S2} implies that the second term on the right hand side 
		is of order $h$, since  
		\begin{align*}
			\langle q_t\sqrt{Q},(P_h-I)B\sqrt{Q}\rangle_{L_2^0} 
			& \leq \|q_t\sqrt{Q}\|_{L_2^0}\|(P_h-I)B\sqrt{Q}\|_{L_2^0} \\
			& \leq Ch \|q_t\sqrt{Q} \|_{L_2^0} \|B\sqrt{Q}\|_{L^0_{2,1}} .
		\end{align*} 
		We thus obtain that  $T_1 \ge - C h$, which together with $T_2\ge 0$, implies the assertion.
	\end{proofThm}

	\subsection{Approximation of the Optimal Feedback Control}
	
	In the whole section we fix $h\in (0,1]$. However, every result of this 
	section remains true in the limit $h\downarrow 0$, i.e. if we consider the 
	optimal control problem \ref{FCP} instead of the finite element approximation 
	\ref{FESCP}. The proofs can be easily adapted to this situation. 
	
	\begin{lemma}\label{ctsLip}
		Under the Assumption \ref{A1}, for any $G_1,G_2\in U_{\mathrm{ad}}^L$ there 
		exists a constant $C>0$ depending only on $T$ and the Lipschitz constant  
		of $G_2$, such that 
		\begin{align*}
			& \E \left[ \sup_{r\in [0,T]} \| u_r^{G_1,h} -u_r^{G_2,h} \|_H^2 
			+ \int_{0}^{T} \| u_r^{G_1,h} - u_r^{G_2,h} \|_{H^1}^2 dr \right] \\
			& \qquad\qquad \leq C \E \left[ \int_0^T \| P_h ( G_1(r,u_r^{G_1,h})  
			- G_2 (r, u_r^{G_1,h})) \|_H^2 dr \right] 
		\end{align*}
		and 
		\begin{align*}
			& \mathbb{E} \left[ \sup_{r\in [0,T]} \| P_h ( G_1 (r, u_r^{G_1,h}) 
			- G_2 (r, u_r^{G_2,h})) \|_H^2 \right] \\
			& \qquad\qquad \leq C \E\left[ \sup_{r\in [0,T]} \| P_h (G_1(r, 
			u_r^{G_1,h})  - G_2 (r, u_r^{G_1,h})) \|_H^2 \right].
		\end{align*}
	\end{lemma}
	
	\begin{proof}
		For any $r\in [0,T]$, the Lipschitz continuity of $G_2$ implies that 
		\begin{equation} 
			\label{ctsLip:eq1}
			\begin{aligned}
				& \| P_h(G_1(r,u_r^{G_1,h})-G_2(r,u_r^{G_2,h}))\|_H\\
				& \leq \| P_h ( G_1(r,u_r^{G_1,h}) 
				-G_2(r,u_r^{G_1,h}))\|_H + \|P_h(G_2(r,u_r^{G_1,h}) 
				-G_2(r,u_r^{G_2,h})) \|_H \\
				& \leq \|P_h(G_1(r,u_r^{G_1,h}) - G_2(r,u_r^{G_1,h})) \|_H 
				+ C \| u_r^{G_1,h}-u_r^{G_2,h} \|_H ,
			\end{aligned}
		\end{equation} 
		where $C$ is the Lipschitz constant of $G_2$. Now by It\^o's formula we get 
		for any $r\in [0,T]$
		\begin{equation}\label{lipF}
			\begin{aligned}
				& \|u_r^{G_1,h}-u_r^{G_2,h}\|_H^2\\
				& =-2\int_{0}^{r}\|u_s^{G_1,h}-u_s^{G_2,h}\|_{H^1}^2 ds \\
				&\quad + 2\int_{0}^{r}\|u_s^{G_1,h}-u_s^{G_2,h}\|_H^2 ds
				+ 2\int_{0}^{r}\langle P_h(\mathcal{F}(u_s^{G_1,h}) 
				- \mathcal{F}(u_s^{G_2,h})),u_s^{G_1,h}-u_s^{G_2,h} \rangle_H ds \\
				& \quad + 2\int_{0}^{t} \langle P_h(G_1(s,u_s^{G_1,h}) 
				- G_2 ( s, u_2^{G_2,h})), u_s^{G_1,h}-u_s^{G_2,h} \rangle_H ds \\
				& \leq -2\int_{0}^{r} \| u_s^{G_1,h} - u_s^{G_2,h}\|_{H^1}^2 ds + 2\int_{0}^{r}\|u_s^{G_1,h}-u_s^{G_2,h}\|_H^2 ds\\
				&\quad + 2 \int_{0}^{r} \| \mathcal{F}(u_s^{G_1,h})  
				- \mathcal{F} (u_s^{G_2,h}))\|_H \|u_s^{G_1,h} - u_s^{G_2,h} \|_H ds \\
				& \quad + 2\int_{0}^{r} \| P_h(G_1(s,u_s^{G_1,h}) 
				- G_2 ( s, u_2^{G_2,h})) \|_H \|u_s^{G_1,h} - u_s^{G_2,h}\|_H ds.
			\end{aligned}
		\end{equation}
		Using Young's inequality and the Lipschitz continuity of $f$ we obtain 
		by the previous considerations
		\begin{align*}
			& \sup_{r\in [0,t]}\|u_r^{G_1,h}-u_r^{G_2,h}\|_H^2  
			+ \int_{0}^{t} \| u_r^{G_1,h} - u_r^{G_2,h}\|_{H^1}^2 dr \\
			& \leq C\int_{0}^t \sup_{r\in [0,s]} \| u_r^{G_1,h}  
			- u_r^{G_2,h} \|_H^2 ds  + \| P_h(G_1(r,u_r^{G_1,h}) - G_2(r,u_r^{G_1,h})) \|_H^2 dr ,  
		\end{align*}
		where the constant $C$ depends on the Lipschitz constants of 
		$G_2$ and of $f$. Gronwalls inequality now implies that 
		\begin{equation} 
			\label{ctsLip:eq2}
			\begin{aligned}
				\sup_{r\in [0,T]} \| u_r^{G_1,h} -u_r^{G_2,h} \|_H^2 
				& + \int_{0}^{T} \| u_r^{G_1,h} - u_r^{G_2,h} \|_{H^1}^2 dr \\
				& \leq C \int_0^T \| P_h ( G_1(r,u_r^{G_1,h})  
				- G_2 (r, u_r^{G_1,h})) \|_H^2 dr .
			\end{aligned} 
		\end{equation} 
		which yields the first inequality taking expectations. 
		
		For the proof of the second inequality we first take the supremum in 
		\eqref{ctsLip:eq1} 
		\begin{align*}
			\sup_{r\in [0,T]} & \| P_h(G_1(r,u_r^{G_1,h})-G_2(r,u_r^{G_2,h}))\|_H \\
			& \leq \sup_{r\in [0,T]} \|P_h(G_1(r,u_r^{G_1,h}) - G_2(r,u_r^{G_1,h})) 
			\|_H 
			+ C \sup_{r\in [0,T]}\| u_r^{G_1,h}-u_r^{G_2,h} \|_H ,
		\end{align*}
		and inserting \eqref{ctsLip:eq2} yields that 
		\begin{align*}
			\sup_{r\in [0,T]} & \| P_h(G_1(r,u_r^{G_1,h})-G_2(r,u_r^{G_2,h}))\|_H \\
			& \leq \sup_{r\in [0,T]} \|P_h(G_1(r,u_r^{G_1,h}) - G_2(r,u_r^{G_1,h})) 
			\|_H 
			\\
			&\quad + C\int_0^T \| P_h ( G_1(r,u_r^{G_1,h})  - G_2 (r, u_r^{G_1,h})) \|_H^2 dr 
			\\
			& \le C \sup_{r\in [0,T]} \|P_h(G_1(r,u_r^{G_1,h}) - G_2(r,u_r^{G_1,h})) 
			\|_H   .
		\end{align*}
		Taking squares and expectations gives the second inequality. 
	\end{proof}
	
	\begin{remark}\label{Rem5}
		The proof of Lemma \ref{ctsLip} can also be easily adapted to the situation where $f$ is only one-sided Lipschitz continuous. In this case the Nemytskii operator $\mathcal{F}$ is also one-sided Lipschitz continuous and since $u^{G_1,h}_s,u^{G_2,h}_s\in S_h$, for all $s\in [0,T]$, we have in \eqref{lipF}
		\begin{align*}
			&\langle P_h(\mathcal{F}(u_s^{G_1,h}) 
			- \mathcal{F}(u_s^{G_2,h})),u_s^{G_1,h}-u_s^{G_2,h} \rangle_H \\
			&=\langle \mathcal{F}(u_s^{G_1,h}) 
			- \mathcal{F}(u_s^{G_2,h}),u_s^{G_1,h}-u_s^{G_2,h} \rangle_H \\
			&\leq C||u_s^{G_1,h}-u^{G_2,h}_s||_H,
		\end{align*}
		where we used the definition of the orthogonal projection $P_h$.
	\end{remark}
	
	\begin{lemma}\label{LipCost2}
		Under the Assumptions \ref{A1}, for any $G_1,G_2\in U_{\mathrm{ad}}^L$ there 
		exists a constant $C>0$ depending only on $T$ and the Lipschitz constants 
		of $G_1$ and $G_2$, such that it holds
		\begin{align*}
			& |J^h(G_1)-J^h(G_2)| \\
			& \leq C\E \left[  \sup_{r\in [0,T]} \|P_h(G_1(r,u_r^{G_1,h}) -
			G_2(r,u_r^{G_1,h})) \|_H^2 \right]^{1/2}.
		\end{align*}
	\end{lemma}
	
	\begin{proof} 
		Lemma \ref{uniformest} implies that  
		\begin{align*} 
			\mathbb{E} \left[  \sup_{t\in [0,T]} \| P_hG_i ( t,u_t^{G_i,h})\|_H^2  +  
			\sup_{t\in [0,T]} \| u_t^{G_i,h} \|_H^2 \right] \leq C,
		\end{align*} 
		for $i=1,2$, for some $C>0$ that depends only on the Lipschitz constants 
		of $G_1$ and $G_2$ and $T$, but is independent of $h$. Therefore a simple 
		application of Lemma \ref{LipCost} yields 
		\begin{align*}
			& |J^h(G_1)-J^h(G_2)|^2\\
			& \leq C \left( \E \left[ \sup_{t\in [0,T]} \|u_t^{G_1,h}  
			- u_t^{G_2,h}\|_H^2\right] 
			+ \mathbb{E} \left[ \sup_{t\in [0,T]}\|P_h(G_1(t,u_t^{G_1,h}) 
			-G_2(t,u_t^{G_2,h})) \|_H^2 \right] \right),
		\end{align*}
		for some constant $C>0$ depending only on $T$ and the Lipschitz constants 
		of $G_1$ and $G_2$. Now Lemma \ref{ctsLip} yields the desired result.
	\end{proof}
	
	We are now in the position to prove our second main result.
	
	\begin{proofThm}{\ref{main2}}
		Let $R > 0$ and
		$$ 
		\mathcal{B}_R := 
		\left\{ \sup_{t \in [0, T]} \|u_t^{\hat{G}, h} \|_H 
		< R\right\} . 
		$$
		Then 
		\begin{align*} 
			\mathbb{P}\left( \mathcal{B}_R^c  \right) 
			& \le \frac{1}{R^2} \E\left(  \sup_{t \in [0, T]} \|u_t^{\hat{G}, h} \|^2_H 
			\right)  \\
			&\le C\frac{1}{R^2},
		\end{align*} 
		for some constant $C$ that is independent of $h$ by Lemma \ref{uniformest}. 
		This implies, again using Lemma \ref{uniformest}, that for any $m\in \mathbb{N}$
		\begin{align*} 
			& \E \left[  \sup_{t\in [0,T]}  \| P_h(G^{h,m}(t,u_t^{\hat{G},h})  
			- \hat{G}(t,u_t^{\hat{G},h}))\|_H^2\right] \\
			& \qquad = \E \left[ \left( \mathbf{1}_{\mathcal{B}_R} 
			+   \mathbf{1}_{\mathcal{B}^c_R} \right) \sup_{t\in [0,T]}   
			\| P_h(G^{h,m}(t,u_t^{\hat{G},h})  
			- \hat{G}(t,u_t^{\hat{G},h}))\|_H^2\right] \\ 
			& \qquad \le \E \left[ \mathbf{1}_{\mathcal{B}_R} 
			\sup_{t\in [0,T]}  \| P_h(G^{h,m}(t,u_t^{\hat{G},h})  
			- \hat{G}(t,u_t^{\hat{G},h}))\|_H^2\right] \\
			& \qquad\qquad  +  \mathbb{P}\left( \mathcal{B}_R^c  \right)^{1/2}   
			\E \left[ \sup_{t\in [0,T]}   
			\| P_h(G^{h,m}(t,u_t^{\hat{G},h})  
			- \hat{G}(t,u_t^{\hat{G},h}))\|_H^4\right]^{1/2} \\  
			&\qquad\le  \epsilon^{h,R}_m+C_h\frac{1}{R},
		\end{align*}
		for some constant $C_h$ independent of $m$. Lemma \ref{LipCost2} now yields
		\begin{align*}
			|J^h(G^{h,m})-J^h(\hat{G})|&\leq \sqrt{\epsilon^{h,R}_m+C_h\frac{1}{R}}.
		\end{align*}
		Using Theorem \ref{FEapprox} we get 
		\begin{align*}
			\inf_{\mathfrak{g} \in \mathbb{A}} J(\mathfrak{g})  
			& \geq J^h(\hat{G}) - C\epsilon_1(h) \\
			& \geq J^h(G^{h,m} ) - C\epsilon_1(h) - \sqrt{\epsilon^{h,R}_m+C_h\frac{1}{R}} \\
			& \geq \inf_{G\in \mathbb{U}^{h,m}}J^h(G) - C\epsilon_1(h) -
			\sqrt{\epsilon^{h,R}_m+C_h\frac{1}{R}},
		\end{align*}
		where 
		$$ 
		\epsilon_1 (h) = C\left( 1+ \E \left[ \int_0^T  \|\hat{G}(t, u_t^{\hat{G}}) 
		\|_{H^1}^2 \, dt \right]^{1/2}\right)h
		$$
		and constants $C,C_h>0$, where 
		the constant $C$ does not depend on $h$. 
	\end{proofThm}
	
	\bigskip 
	The lower bound is a simple consequence of Theorem 
	\ref{convexEst}.
	
	\begin{proofThm}{\ref{main3}}
		By Theorem \ref{convexEst} we get 
		\begin{align*}
			\inf_{\mathfrak{g} \in \mathbb{A}} J(\mathfrak{g}) 
			& \leq \inf_{G\in U_{\mathrm{ad}}}J^h(G) + Ch
		\end{align*}
		for some constant $C$ independent of $h$, and since 
		$\mathbb{U}\subseteq U_{\mathrm{ad}}$ we obtain 
		\begin{align*}
			\inf_{\mathfrak{g} \in \mathbb{A}} J(\mathfrak{g}) 
			& \leq \inf_{G\in \mathbb{U}}J^h(G) + Ch.
		\end{align*}
	\end{proofThm}

	\section{Construction of Ansatz Spaces}\label{ansatz}
	
	In this section we will provide a method to construct ansatz spaces $\mathbb{U}$, $\mathbb{U}^{h,m}$ that satisfy the uniform or uniform Lipschitz approximation property with respect to an optimal control. Although optimal controls are typically not known explicitly, it is still possible to construct explicit ansatz spaces if the optimal feedback $\hat{G}$ is assumed to be continuous or Lipschitz continuous. The core idea is to consider a suitable dense subset of $\mathcal{C}([0,T]\times \mathbb{R}^{N_h},\mathbb{R}^{N_h})$ with respect to the topology of compact convergence, in order to approximate the continuous functions
	\begin{align*}
		g^h_i(t,u_1,...,u_{N_h})=\langle \hat{G}(t,\sum_{j=1}^{N_h}u_je_j),e_i\rangle_H,
	\end{align*} 
	where $N_h$ denotes the dimension of some finite dimensional subspace $S_h\subseteq H^1$ with orthonormal basis $\{e_1,...,e_{N_h}\}$. Based on this approximation, the ansatz spaces $\mathbb{U}$, $\mathbb{U}^{h,m}$ can be explicitly constructed. Examples for dense subsets suitable for the numerical implementation will be given in Section \ref{ex}. 
	
	In the following let $(S_h)_{h\in (0,1]}$ be a family of finite dimensional subspaces of $H^1$ with orthogonal projections $P_h$ satisfying 
	\begin{align*}
		||P_hu-u||_H\rightarrow 0,\text{ as }h\rightarrow 0.
	\end{align*}
	By $N_h$ we denote the dimension of $S_h$ and by $\{e_1,...,e_{N_h}\}$ we denote an orthonormal basis of $S_h$. Furthermore for any $h\in (0,1]$ we define the function
	\begin{align*}
		g^h&:[0,T]\times \mathbb{R}^{N_h}\rightarrow \mathbb{R}^{N_h}\\
		g_i^h(t,u)&:=\langle \hat{G}(t,\sum_{j=1}^{N_h}u_je_j),e_i\rangle_H, \quad i=1,...,N_h.
	\end{align*} 
	
	\subsection{Uniform Ansatz Space $\mathbb{U}$}\label{uniformAnsatz}
	
	In the first part we consider the case of bounded controls, i.e. the  control space is given by $\mathcal{U}=\mathcal{B}_H(0,K)$, for some $K>0$. In this situation we will construct an ansatz space satisfying the uniform approximation property. 
	
	Since $\hat{G}$ is assumed to be continuous, the functions $(g^h)_{h\in (0,1]}$ are also continuous. In particular, it is possible to approximate these functions by simpler functions that can be treated numerically, e.g. artificial neural networks. In the following we consider for any $h\in (0,1]$ a set of Lipschitz approximations $\mathcal{N}^{h}\subseteq \mathcal{C}^{1-\text{Höl}}([0,T]\times \mathbb{R}^{N_h},\mathbb{R}^{N_h})$, such that there exists a sequence $(\psi^{h,k})_{k\in \mathbb{N}}\subseteq \mathcal{N}^h$ with 
	\begin{align*}
		\lim\limits_{k\rightarrow \infty}\sup_{(t,u)\in [0,T]\times \mathcal{B}_{\mathbb{R}^{N_h}}(0,R)}|\psi^{h,k}(t,u)-g^h(t,u)|=0,
	\end{align*} 
	for all $R>0$. For a particular choice of $\mathcal{N}^{h}$ we refer to our examples in Section \ref{ex}. Now we define the ansatz space
	\begin{multline*}
		\mathbb{U}:=\bigg \{ G\left( t,\sum_{i=1}^\infty u_i e_i\right) (x) 
		= \sum_{i=1}^{N_h}C\psi_i(t,\eta^l(u_1,\dots,u_n)) e_i(x)
		\bigg |\\
		\sup_{t\in [0,T]}\sup_{u\in H}\|G(t,u)\|_H\leq K, \text{ where }\psi\in \mathcal{N}^{h},h\in (0,1],l 
		\in \mathbb{N}, C\in \mathbb{R} \bigg \},
	\end{multline*}
	where $\psi_i$ denotes the $i$-th component function of $\psi$ and for $l\in \mathbb{N}$ the function $\eta^l : \mathbb{R}^{N_h}\rightarrow \mathbb{R}^{N_h}$ 
	\begin{align*}
		\eta^l (x) 
		& = \begin{cases}
			x & |x|\leq l \\
			l\frac{x}{|x|} & |x|>l
		\end{cases}
	\end{align*}
	is a cutoff function. 
	
	Since the elements of $\mathbb{U}$ are Lipschitz continuous, equation \eqref{spdeFeedback2} has a unique strong solution for every $G\in \mathbb{U}$, hence $\mathbb{U}\subseteq U_{ad}$. In order to show that the constructed ansatz space satisfies the uniform approximation property with respect to $\hat{G}$, we need to construct a family $(G^{h,m})_{h\in (0,1],m\in \mathbb{N}}$ that satisfies \eqref{u1} and \eqref{u2}. Recalling the assumption on $\mathcal{N}^{h}$, we can find a sequence $(\psi^{h,k_m})_{m\in \mathbb{N}} \subseteq \mathcal{N}^{h}$ with 
	\begin{align}
		\sup_{(t,u)\in [0,T]\times \mathcal{B}_{\mathbb{R}^{N_h}}(0,m)}|\psi^{h,k_m}(t,u)-g^h(t,u)|\leq \frac{K}{m}.
	\end{align}
	Now we define the family $(G^{h,m})_{h\in (0,1],m\in \mathbb{N}}$ by
	\begin{equation} 
		G^{h,m}(t,u) 
		:= \left( 1-\frac{1}{m} \right) 
		\sum_{i=1}^{N_h} \psi^{h,k_m}_i (t,
		\eta^m(\langle u,e_1\rangle_H ,\ldots, \langle u,e_{N_h}\rangle_H)) e_i.
	\end{equation}
	It clearly holds $\psi^{h,k_m}(t,\eta^m(u))=\psi^{h,k_m}(t,u)$ on  
	$[0,T]\times \mathcal{B}_{\mathbb{R}^{N_h}}(0,m)$ and for any  
	$(t,u)\in[0,T]\times\mathcal{B}_{\mathbb{R}^{N_h}}(0,m)$
	\begin{align*}
		|\psi^{h,k_m}(t,u)| 
		& \leq |g^h (t,u)| + |\psi^{h,k_m}(t,u) - g^h (t,u)| \\
		& \leq K + \frac{K}{m} 
		\leq \left( 1 + \frac 1m\right)K , 
	\end{align*}
	since 
	$$ 
	|g^h (t,u)| \le \|P_{h} \hat{G} (t, \sum_{j=1}^{N_h} 
	u_j e_j ) \|_H \le \| \hat{G} (t, \sum_{j=1}^{N_h}  
	u_j e_j ) \|_H\le K.
	$$
	Furthermore, on $[0,T]\times\mathcal{B}_{\mathbb{R}^{N_h}}(0,m)^c$ we have 
	\begin{align*}
		|\psi^{h,k_m}(t,\eta^m(u))| 
		& = |\psi^{h,k_m}(t,m\frac{u}{|u|})| \\
		& \leq |g^h(t,m\frac{u}{|u|})|+|\psi^{h,k_m}(t,m\frac{u}{|u|}) 
		- g^h(t,m\frac{u}{|u|})| 
		\\
		& \leq \left( 1+ \frac 1m\right) K.
	\end{align*}
	Therefore it holds
	\begin{align*}
		\| G^{h,m}(t,u) \| 
		& = \left( 1 - \frac{1}{m}\right) |\psi^{h,k_m} 
		\left( t, \eta^m(\langle u,e_1\rangle_H , \ldots , \langle u,e_{N_n}  
		\rangle_H)\right) | \\ 
		& \le \left( 1-\frac 1m\right)\left( 1 + \frac 1m\right) K \le K, 
	\end{align*}
	hence $G^{h,m}\in \mathbb{U}$. Now for any $R>0$ and any $\epsilon>0$ there exists an $M\in \mathbb{N}$, such 
	that $\mathcal{B}_H (0,R)\subset  \mathcal{B}_H (0,m)$ and 
	$\frac{2K}{m}<\epsilon$ for every $m\geq M$. Therefore, since 
	\begin{align*}
		\langle \hat{G}^h(t,u),e_i\rangle_H&=\langle P_h\hat{G}(t,P_hu),e_i\rangle_H\\
		&=g_i^h(t,\langle u,e_1\rangle_H,...,\langle u,e_{N_h}\rangle_H),
	\end{align*} 
	for $i=1,...,N_h$, we have for any $m\geq M$
	\begin{align*}
		& \sup_{(t,u)\in [0,T]\times \mathcal{B}_H (0,R)} \| G^{h,m}(t,u) 
		- \hat{G}^h(t,u) \| \\
		& \leq \sup_{(t,u)\in [0,T]\times \mathcal{B}_H (0,m)} 
		\|G^{h,m}(t,u) - \hat{G}^h(t,u) \| \\
		& \leq \sup_{(t,u)\in [0,T] \times \mathcal{B}_{\mathbb{R}^{N_h}}(0,m)} 
		|\left( 1- \frac 1m\right) \psi^{h,k_m}(t,u) - g^h(t,u)| \\
		& \leq  \sup_{(t,u)\in [0,T]\times \mathcal{B}_{\mathbb{R}^{N_h}}(0,m)} 
		\left( 1- \frac 1m\right)  | \psi^{h,k_m}(t,u)-g^h(t,u)|  
		+ \frac 1m |g^h (t,u)| \\
		& < \frac{2K}{m} < \epsilon.
	\end{align*}

	\subsection{Ansatz Spaces of Fixed Dimension $\mathbb{U}^{h,m}$}\label{ansatzFix}
	
	In the second part we focus on the case of unbounded controls, i.e. the control space is given by $\mathcal{U}=H$. We assume, that there exists a unique 
	optimal control given in feedback form by
	\begin{align*}
		\mathfrak{\hat{g}}_t=\hat{G}(t,u_t^{\hat{G}}),
	\end{align*} 
	for some feedback $\hat{G}:[0,T]\times H\rightarrow H$ which is 
	Lipschitz continuous in $u$ and where $u^{\hat{G}}$ is the unique strong 
	solution to equation (\ref{spdeFeedback}). 
	
	The existence of a Lipschitz continuous feedback control can be ensured, if the solution to the HJB equation (\ref{HJB}) has bounded second derivatives on $[0,T-\epsilon]\times H$, for all $\epsilon\in\, ]0,T[$, see \cite[Theorem~4.201]{FGS17a}. Sufficient conditions for this assumption to be true are given in \cite[Theorem~4.155]{FGS17a}. 
	
	In the above setting we will construct a sequence of ansatz spaces that satisfies a uniform Lipschitz approximation property of order $h\in (0,1]$ with respect to $\hat{G}$. Therefore we consider Lipschitz approximations $(\mathcal{N}^{h,k})_{k\in \mathbb{N}}$, $\mathcal{N}^{h,k}\subseteq \mathcal{C}^{1-\text{Höl}}([0,T]\times \mathbb{R}^{N_h},\mathbb{R}^{N_h})$, such that there exist a sequence of Lipschitz continuous functions $\psi^{h,k}\in \mathcal{N}^{h,k}, k\in \mathbb{N}$ with Lipschitz constant independent of $k$, such that for all $R>0$ 
	\begin{align*}
		\lim\limits_{k\rightarrow \infty}\sup_{(t,u)\in [0,T]\times \mathcal{B}_{\mathbb{R}^{N_h}}(0,R)}|\psi^{h,k}(t,u)-g^h(t,u)|= 0.
	\end{align*}
	Then we define the ansatz space
	\begin{align*}
		\mathbb{U}^{h,k}:=\bigg \{ G\left( t,\sum_{i=1}^\infty u_i e_i\right) (x) 
		= \sum_{i=1}^{N_h}\psi_i(t,(u_1,\dots,u_n)) e_i(x)
		\bigg |
		\psi\in \mathcal{N}^{h,k}\bigg \}.
	\end{align*}
	Again, since the elements of $\mathbb{U}^{h,k}$ are Lipschitz continuous, equation \eqref{spdeFeedback2} has a unique strong solution for every $G\in \mathbb{U}^{h,k}$, in particular $\mathbb{U}^{h,k}\subseteq U_{ad}^L$. By defining $G^{h,m}$ as 
	\begin{equation} 
		G^{h,k}(t,u) 
		:=
		\sum_{i=1}^{N_h} \psi^{h,k}_i (t,
		(\langle u,e_1\rangle_H ,\ldots, \langle u,e_{N_h}\rangle_H)) e_i,
	\end{equation}
	it is not difficult to verify, that $(\mathbb{U}^{h,k})_{k\in \mathbb{N}}$ satisfies a uniform Lipschitz approximation property with respect to $\hat{G}$, if we can show that $G^{h,m}$ is Lipschitz continuous with Lipschitz constant independent of $m$.
	
	By the assumption on $\mathcal{N}^{h,k}$ we have for any $s,t\in [0,T]$ and $u,v\in H$ 
	\begin{align*}
		&||G^{h,m}(t,u)-G^{h,m}(s,v)||_H\\
		&=|\psi^{h,m}(t,\langle u,e_1\rangle_H,...,\langle u,e_{N_h}\rangle_H)-\psi^{h,m}(s,\langle v,e_1\rangle_H,...,\langle v,e_{N_h}\rangle_H)|\\
		&\leq L_h\left(|t-s|+|(\langle u-v,e_1\rangle_H,...,\langle u-v,e_{N_h}\rangle_H)|\right)\\
		&\leq L_h\left(|t-s|+||u-v||_H\right),
	\end{align*}
	for some $L_h>0$ independent of $m$. This shows that $G^{h,m}$ is Lipschitz continuous with Lipschitz constant independent of $m$ and therefore $(\mathbb{U}^{h,k})_{k\in \mathbb{N}}$ satisfies a uniform Lipschitz approximation property with respect to $\hat{G}$.
	
	\begin{remark}
		This allows in particular to construct an ansatz space that satisfies the uniform approximation property with respect to $\hat{G}$ in the case of unbounded controls by
		\begin{align*}
			\mathbb{U}:=\bigcup_{h\in (0,1]}\bigcup_{k=1}^\infty \mathbb{U}^{h,k}
		\end{align*}
	\end{remark}
	
	\section{Examples for Ansatz Spaces}
	\label{ex}
	
	In this section we will provide some explicit examples for ansatz spaces to demonstrate the scope of our two main theorems.
	
	In the first example regarding Theorem \ref{main1} we will give an example for an ansatz space of approximating controls $\mathbb{U}$ that satisfies the uniform approximation property with respect to the optimal feedback $\hat{G}$ by using artificial neural networks with one layer. In the second example regarding Theorem \ref{main2} we will give an example for a sequence of ansatz spaces $(\mathbb{U}^{h,m})_{m\in \mathbb{N}}$ that satisfies the uniform Lipschitz approximation property with respect to the optimal feedback $\hat{G}$ using Gaussian radial basis neural networks and specify explicit approximation rates.
	
	\subsection{Universal Approximation by Neural Networks}
	
	In the first example we will use artificial neural networks to define ansatz spaces that satisfy the uniform and uniform Lipschitz approximation property respectively. Regarding the uniform approximation property we consider the ansatz space $\mathbb{U}$ constructed in Subsection \ref{ansatz} for the approximating set 
	\begin{align*}
		\mathcal{N}^{h}:=\bigcup_{k=1}^{\infty}\mathcal{N}_{k}^{N_h},
	\end{align*}
	where 
	\begin{align*}
		\mathcal{N}_k^n:=\big \{\psi(x)=C \theta (A x+b)
		\big |\,A\in \mathbb{R}^{k\times n},b\in \mathbb{R}^k,C \in 
		\mathbb{R}^{n\times k} \big \}
	\end{align*}
	denotes the set of all 1-layer artificial neural networks from 
	$\mathbb{R}^n$ to $\mathbb{R}^n$ with $k$ neurons, for a given 
	non-polynomial, Lipschitz continuous activator function $\theta$, where the activator function is evaluated componentwise. Thanks to the discussion in Section \ref{uniformAnsatz}, we just need to show that for any $h\in (0,1]$ there exists a sequence $(\psi^{h,k})_{k\in \mathbb{N}}\subseteq \mathcal{N}^h$ with 
	\begin{align*}
		\lim\limits_{k\rightarrow \infty}\sup_{(t,u)\in [0,T]\times \mathcal{B}_{\mathbb{R}^{N_h}}(0,R)}|\psi^{h,k}(t,u)-g^h(t,u)|=0,
	\end{align*} 
	for all $R>0$. 
	
	Since for every $h\in (0,1]$ the function $g^h$ is continuous, this however is a simple consequence from the classical universal approximation 
	result \cite[Theorem~3.1]{Pin99}.\\
	
	For the uniform Lipschitz approximation property we consider for fixed $h\in (0,1)$ the ansatz spaces $(\mathbb{U}^{h,k})_{k\in \mathbb{N}}$ constructed in Subsection \ref{ansatzFix} for the approximating sets $(\mathcal{N}^{N_h}_k)_{k\in \mathbb{N}}$. By \cite[Proposition~10]{CL22} there exists a sequence $(\psi^{h,k})_{k\in \mathbb{N}}$ of artificial neural networks  with Lipschitz constant independent of $k$ and $\psi^{h,k}\in \mathcal{N}^{h,k}$, such that 
	\begin{align*}
		\sup_{(t,u)\in [0,T]\times \mathcal{B}_{\mathbb{R}^{N_h}}(0,k^{\frac{1}{3(N_h+1)}})}|\psi^{h,k}(t,u)-g^h(t,u)|\leq k^{-\frac{1}{3(N_h+1)}}\rightarrow 0,
	\end{align*}
	as $k\rightarrow \infty$. This now implies that the corresponding ansatz spaces $(\mathbb{U}^{h,k})_{k\in \mathbb{N}}$  satisfy the uniform Lipschitz approximation property.
	
	\subsection{Interpolation by Radial Basis Functions}
	
	In our second example we will interpolate the finitely based approximations 
	by Gaussian radial basis function neural networks and define a suitable 
	ansatz space of controls for the approximation that satisfies a uniform 
	Lipschitz approximation property with respect to the optimal feedback 
	$\hat{G}$. The main idea is to determine the optimal feedback at some suitable 
	states $(t_1,x_1),...,(t_N,x_N)$ and interpolate afterwards. In the following we introduce the notation from \cite{Wen04}. For some 
	$\kappa > 0$ let 
	\begin{align*}
		\Phi(x,y) := \exp(-\kappa |x-y|^2) 
	\end{align*}
	denote the Gaussian kernel. Furthermore we define for $O\subseteq 
	\mathbb{R}^{d}$ 
	\begin{align*}
		F_\Phi (O) = \text{span} \{ \Phi (\cdot , y) 
		\mid y\in O \} , 
	\end{align*}
	equipped with the scalar product 
	\begin{align*}
		\langle \sum_{i=1}^{N} \alpha_i \Phi(\cdot,x_j), \sum_{j=1}^{M} \beta_j 
		\Phi(\cdot,y_j)\rangle_\Phi  
		:= \sum_{i=1}^{N}\sum_{j=1}^{M}\alpha_i\beta_j\Phi(x_i,y_j).
	\end{align*}
	We denote by $\mathcal{F}_{\Phi}(O)$ the completion of $F_{\Phi}(O)$ with respect to the norm $\|\cdot\|_{\Phi}$ induced by $\langle \cdot,\cdot\rangle_{\Phi}$. Now we define the 'point evaluation' map 
	\begin{align*}
		\mathcal{R}: \mathcal{F}_\Phi( O)\rightarrow \mathcal{C} ( O), 
		\mathcal{R}(f)(x) := \langle f, \Phi(\cdot,x)\rangle_\Phi
	\end{align*}
	to identify abstract elements of $\mathcal{F}_\Phi(O)$ with functions 
	and introduce the corresponding native Hilbert space of $\Phi$ by 
	\begin{align*}
		\mathcal{N}_\Phi(O):=\mathcal{R}(\mathcal{F}_\Phi(O)),
	\end{align*}
	with the inner product
	\begin{align*}
		\langle f_1,f_2\rangle_{\mathcal{N}_\Phi(O)} 
		:= \langle \mathcal{R}^{-1}f_1,\mathcal{R}^{-1}f_2\rangle_\Phi.
	\end{align*}
	If $O=\mathbb{R}^{d}$, the native space is given by 
	\begin{align*}
		\mathcal{N}_\Phi(\mathbb{R}^{d}) 
		= \{ f\in L^2(\mathbb{R}^{d}) \cap \mathcal{C}(\mathbb{R}^{d}) 
		\mid \hat{f}/\hat{\Phi}^{1/2} \in L^2(\mathbb{R}^{d})\},
	\end{align*}
	with the inner product
	\begin{align*}
		\langle f_1,f_2\rangle_{\mathcal{N}_\Phi (\mathbb{R}^{d})} 
		= (2\pi)^{-d/2} \langle \hat{f}_1 / \hat{\Phi}^{1/2}, 
		\hat{f}_2/ \hat{\Phi}^{1/2} \rangle_{L^2(\mathbb{R}^{d})} , 
	\end{align*}
	where $\hat{f}$ denotes the analytic Fourier transform of $f$. 
	
	For a given function $g\in \mathcal{N}_\Phi(O)$ we will consider approximations of the type
	\begin{equation}\label{interpolant}
		s_{g,X}(x):=\sum_{k=1}^{K}\alpha_k \Phi(x,x_k),
	\end{equation}
	for a discrete set $X=\{x_1,...,x_K\}\subseteq O$ and some $\alpha_1,...,\alpha_K\in \mathbb{R}$, such that 
	\begin{equation}\label{interpolProp}
		s_{g,X}(x_k)=g(x_k),
	\end{equation}
	for $k=1,...,K$. 
	
	\begin{lemma}\label{Lipconstraidal}
		Let $g\in \mathcal{N}_\Phi(\mathbb{R}^{d})$, $X=\{x_1,...,x_K\} 
		\subseteq O\subseteq \mathbb{R}^{d}$ and $s_{g,X}$ be of the form 
		\eqref{interpolant}, satisfying \eqref{interpolProp} for  
		$g|_{O}\in \mathcal{N}_\Phi(O)$. Then $s_{g,X}$ is Lipschitz continuous with Lipschitz constant given by $2\kappa \| g \|_{\mathcal{N}_\Phi(\mathbb{R}^{d})}^2$.
	\end{lemma}
	
	\begin{proof}
		We first observe that for any $x\in \mathbb{R}^{d}$
		\begin{align*}
			s_{g,X}(x) 
			=\langle s_{g,X},\Phi(\cdot,x)\rangle_{\Phi} 
		\end{align*}
		and therefore for $x,y\in \mathbb{R}^d$
		\begin{align*}
			|s_{g,X}(x)-s_{g,X}(y)|
			& =|\langle s_{g,X}, \Phi (\cdot,x)-\Phi (\cdot,y)\rangle 
			_{\Phi}|.
		\end{align*}
		Using the Cauchy-Schwarz inequality we get
		\begin{align*}
			|s_{g,X}(x)-s_{g,X}(y)|^2&\leq \| s_{g,X} \|_\Phi^2\| 
			\Phi(\cdot,x)-\Phi(\cdot,y)\|_{\Phi}^2.
		\end{align*}
		Now 
		\begin{align*}
			\| 
			\Phi(\cdot,x)-\Phi(\cdot,y)\|_{\Phi}^2
			&\leq 4\kappa^2|x-y|^2.
		\end{align*}
		Furthermore we have 
		\begin{align*}
			\| s_{g,X} \|_\Phi^2 
			&\leq \| g \|_{\mathcal{N}_\Phi(\mathbb{R}^{d})}^2.
		\end{align*}
		This concludes the proof of this lemma. 
	\end{proof}
	
	In the following we define for $X = \{ x_1 , \ldots , x_K\} 
	\subseteq 
	O$
	\begin{align*}
		h_{X,O}:=\sup_{x\in O}\min_{1\leq j\leq K}|x-x_j|.
	\end{align*}
	
	Then we consider the following two slightly modified results from \cite{Wen04}. The first result is an immediate consequence of the proof of \cite[Proposition~14.1]{Wen04}.
	
	\begin{lemma}\label{grid}
		Let $O\subseteq \mathcal{B}_R(x_0)$, for some $R>0$ and  
		$X=\{x_1,...,x_K\}\subseteq O$ be quasi uniform (q.u.) with respect 
		to $c_{qu}>0$, i.e. 
		\begin{align*}
			q_X \leq h_{X,O}\leq c_{qu}q_X,
		\end{align*}
		where $q_X:=\frac{1}{2}\min_{i\not=j}|x_i-x_j|$. Then it holds 
		\begin{align*}
			h_{X,O}\leq 2Rc_{qu}K^{-1/d}.
		\end{align*}
	\end{lemma}
	
	The second result is a consequence of \cite[Proposition~11.14]{Wen04}
	
	\begin{theorem}\label{InterpolEst}
		Let $f\in \mathcal{N}_\Phi(O)$ and $s_{f,X}$ denote its interpolant 
		based on the quasi uniform set $X=\{x_1,...,x_K\}\subseteq O$.
		Then there exists a constant $C > 0$, such that 
		\begin{align*}
			\sup_{x\in O}|f(x)-s_{f,X}(x)| 
			\leq C^lh_{X,O}^l \|f\|_{\mathcal{N}_\Phi(O)},
		\end{align*}
		for any $l\in \mathbb{N}$.
	\end{theorem}
	
	Now we consider for fixed $h\in (0,1]$ and $c_{qu}>0$, the sequence of approximations 
	
	\begin{multline*}
		\mathcal{N}^{h,k} := \{ s(t,x) = \sum_{i=1}^k \alpha_i \Phi((t,x),  
		( t_i,x_i))\\
		|((t_1,x_1),...,(t_k,x_k))\in  [0,T]\times 
		\mathcal{B}_{\mathbb{R}^{N_h}}(0,k^{1/(2N_h)})\text{ q.u.}, \alpha_1,...,\alpha_k\in \mathbb{R}\}.
	\end{multline*}
	
	Thanks to the discussion in Subsection \ref{ansatzFix}, we just need to show, that  there exist a sequence of Lipschitz continuous functions $s^{h,k}\in \mathcal{N}^{h,k}, k\in \mathbb{N}$ with Lipschitz constant independent of $k$, such that for all $R>0$ 
	\begin{align*}
		\lim\limits_{k\rightarrow \infty}\sup_{(t,u)\in [0,T]\times \mathcal{B}_{\mathbb{R}^{N_h}}(0,R)}|s^{h,k}(t,u)-g^h(t,u)|= 0.
	\end{align*}
	
	In the following impose stronger assumptions on the optimal feedback $\hat{G}$, in particular we assume that the functions $g^h_i$, $i=1,\dots N_h$ are elements of $\mathcal{N}_\Phi([0,T]\times\mathbb{R}^{N_h})$. 
	
	Let $l\in \mathbb{N}$, then for any $R>0$ there exists an $M\in \mathbb{N}$, such that 
	$\mathcal{B}_{\mathbb{R}^{N_n}}(0,R)\subseteq  
	\mathcal{B}_{\mathbb{R}^{N_n}}(0,k^{1/(2N_h)})$ for all $k\geq M$.  
	Furthermore by Theorem \ref{InterpolEst} and Lemma \ref{grid} there exists 
	for all $k\geq M$ an element $s^{h,k}\in \mathcal{N}^{h,k}$, such that 
	\begin{align*}
		&\sup_{(t,u)\in [0,T]\times \mathcal{B}_{\mathbb{R}^{N_h}} 
			(0,R)}|s^{h,k}(t,u) - g^h(t,u)| \\
		& \leq \sup_{(t,u)\in [0,T]\times \mathcal{B}_{\mathbb{R}^{N_h}}
			(0,k^{1/(2N_h)})}|s^{h,k}(t,u) - g^h(t,u)|\\
		& \le C^l c_{qu}^l k^{ - l/(2N_h)}\rightarrow 0,
	\end{align*}
	as $k\rightarrow \infty$. Due to Lemma \ref{Lipconstraidal} any $s^{h,k}$ has 
	Lipschitz constant independent of $k$. All together the sequence of ansatz spaces $(\mathbb{U}^{h,k})_{k\in \mathbb{N}}$ satisfies the uniform Lipschitz approximation property of order $h$ with respect to $\hat{G}$.

	\section{Numerical Example}\label{numerics}
	
	In this section we consider the controlled stochastic heat equation in order to validate our algorithm by comparing with the optimal feedback control obtained from the associated Riccati equation, see equation \eqref{Ricc}. The controlled state equation is given by
	\begin{equation}\label{A:stateequation}
		\begin{cases}
			\mathrm{d}u^{\mathfrak{g}}_t = [\Delta u^{\mathfrak{g}}_t + \mathfrak{g}_t]  \mathrm{d}t + 0.01B\mathrm{d}W_t,\quad t\in [0,20]\\
			u^{\mathfrak{g}}_0=u\in L^2(0,20),
		\end{cases}
	\end{equation}
	with Neumann boundary conditions and $u = \mathbf{1}_{[20/3,40/3]}$, where $W$ is a cylindrical Wiener process on $H$ with covariance operator $Q=Id_H$. Here the Hilbert Schmidt operator $B\in L_2(L^2(\Lambda),L^2(\Lambda))$ is given by 
	\begin{align*}
		B=\Delta^{-\gamma},
	\end{align*}
	for $\gamma=0.751$, hence $\|B \sqrt{Q}\|_{L_{2,1}^0}<\infty$. We consider the problem of steering the solution of the stochastic heat equation into the constant zero profile. To this end, we introduce the cost functional
	\begin{equation}
		J(\mathfrak{g}) = \frac{1}{2} \mathbb{E} \left [ \int_0^{20} \| u^{\mathfrak{g}}_t\|_{L^2(0,20)}^2 + \| \mathfrak{g}_t \|_{L^2(0,20)}^2 \mathrm{d}t \right ].
	\end{equation} 
	Note that the second term is a regularization, which is necessary in linear quadratic control theory. We approximate the Riccati equation \eqref{Ricc} numerically, based on $N_h=400$ Fourier coefficients by 
	\begin{equation}\label{RiccApprox}
		\begin{cases}
			\partial_t P^h(t) + P^h(t) \Delta_h +\Delta_h P^h(t) - I_h + (P^h)^2(t) = 0,\;\; t\in [0,T]\\
			P^h(T) = - I_h,
		\end{cases}
	\end{equation}
	to obtain the approximated optimal cost 
	\begin{align*}
		J_{\text{opt}}\approx J(\mathfrak{g}^{Ric,h})= 5.34
	\end{align*} 
	for the feedback control 
	\begin{equation*}
		\mathfrak{g}^{\text{Ric,h}}_t = P^h(t) u^{\mathfrak{g}^{\text{Ric,h}}}_t,
	\end{equation*}
	where $P^h$ is the solution to the associated approximated Riccati equation, and use this approximation as a benchmark.
	
	For the approximation of the optimal control, we use the ansatz space constructed in Subsection \ref{ansatzFix} with respect to the approximating sets of 1-layer artificial neural networks 
	\begin{align*}
		\mathcal{N}^{h,k}:=\bigg \{\psi(t,u)=C\theta\left (A \begin{pmatrix}t\\ u
		\end{pmatrix}+a\right )
		\bigg |\,A\in \mathbb{R}^{k\times (N_h+1)},C\in \mathbb{R}^{(N_h+1)\times k},a\in \mathbb{R}^{k} \bigg \},
	\end{align*}
	in dimension $N_h=400$ with $k=400$ neurons, and ReLU activator function $\theta$. 
	
\begin{figure}[H]
		\begin{minipage}[t]{.49\linewidth} 
			\centering
			\includegraphics[width=\linewidth]{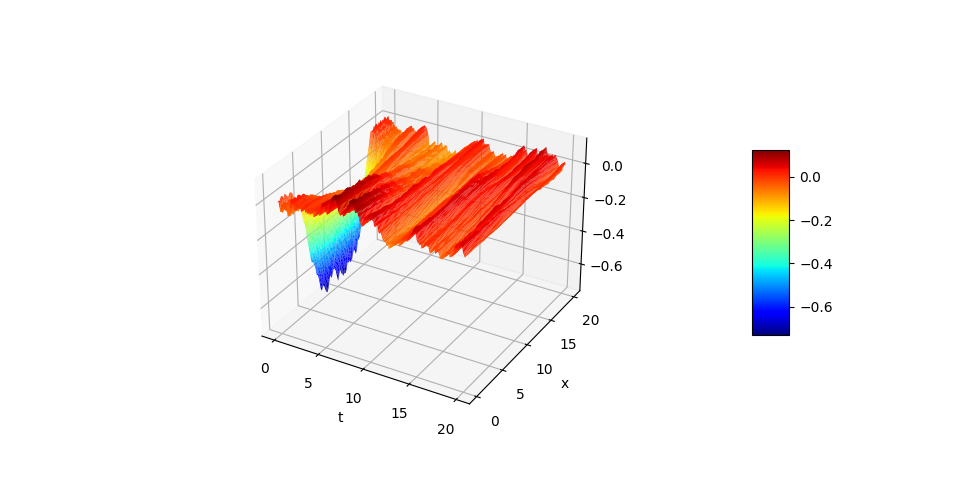}
			(a) Sample approximated control $\mathfrak{g}^{approx}$
		\end{minipage}
		\begin{minipage}[t]{.49\linewidth} 
			\centering
			\includegraphics[width=\linewidth]{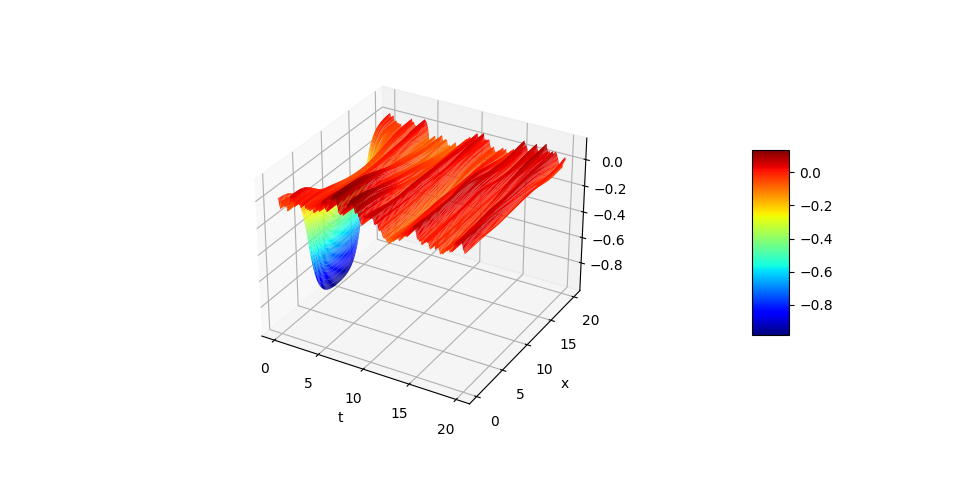}
			(b) Sample optimal control $\mathfrak{g}^{Ric,h}$
		\end{minipage}
		\caption{LQ optimal control}\label{LQFig}
	\end{figure}

	The following shows a realization of the approximated optimal control
	\begin{align*}
		\mathfrak{g}^{\text{approx}}_t=G^{\text{approx}}(t,u_{t}^{G^{\text{approx}}}),
	\end{align*} 
	for $G^{\text{approx}}\in \mathbb{U}^{h,k}$, see Figure \ref{LQFig} (a), and a realization of $\mathfrak{g}^{\text{Ric},h}$ with respect to the same noise realization, see Figure \ref{LQFig} (b).
	
	We ended up with an approximated optimal cost of 
	\begin{align*}
		J_{approx}=J(\mathfrak{g}^{approx})\approx 5.43.
	\end{align*}
	
	\bibliography{spde4}
	\bibliographystyle{abbrv}
	
\end{document}